\documentclass[12pt]{article}
\usepackage{microtype} \DisableLigatures{encoding = *, family = * }
\usepackage{graphicx,subeqn,multirow,slashbox}
\usepackage{ifpdf}
\ifpdf  \fi
\usepackage{amsfonts,amssymb}
\usepackage{natbib}
\bibpunct{(}{)}{;}{a}{,}{,}
\usepackage[bookmarksnumbered=true, pdfauthor={Wen-Long Jin}]{hyperref}

\newcommand{\commentout}[1]{}

\newcommand{\ba}{\begin{array}}
        \newcommand{\ea}{\end{array}}
\newcommand{\bc}{\begin{center}}
        \newcommand{\ec}{\end{center}}
\newcommand{\bdm}{\begin{displaymath}}
        \newcommand{\edm}{\end{displaymath}}
\newcommand{\bds} {\begin{description}}
        \newcommand{\eds} {\end{description}}
\newcommand{\ben}{\begin{enumerate}}
        \newcommand{\een}{\end{enumerate}}
\newcommand{\beq}{\begin{equation}}
        \newcommand{\eeq}{\end{equation}}
\newcommand{\bfg} {\begin{figure}}
        \newcommand{\efg} {\end{figure}}
\newcommand{\bi} {\begin {itemize}}
        \newcommand{\ei} {\end {itemize}}
\newcommand{\bpp}{\begin{pspicture}}
        \newcommand{\epp}{\end{pspicture}}
\newcommand{\bqn}{\begin{eqnarray}}
        \newcommand{\eqn}{\end{eqnarray}}
\newcommand{\bqs}{\begin{eqnarray*}}
        \newcommand{\eqs}{\end{eqnarray*}}
\newcommand{\bsq}{\begin{subequations}}
        \newcommand{\esq}{\end{subequations}}
\newcommand{\bsl} {\begin{slide}[8.8in,6.7in]}
        \newcommand{\esl} {\end{slide}}
\newcommand{\bss} {\begin{slide*}[9.3in,6.7in]}
        \newcommand{\ess} {\end{slide*}}
\newcommand{\btb} {\begin {table}[h]}
        \newcommand{\etb} {\end {table}}
\newcommand{\bvb}{\begin{verbatim}}
        \newcommand{\evb}{\end{verbatim}}

\newcommand {\pd}[2] {{\frac {\partial {#1}} {\partial {#2}}}}

\newcommand{\cas}[1]{{{\left \{ \ba #1 \ea \right. }}}

\newcommand{\reff}[1] {{{Figure \ref {#1}}}}
\newcommand{\refe}[1] {{{(\ref {#1})}}}

\newtheorem{theorem}{Theorem}[section]

\def\pmb#1{\setbox0=\hbox{$#1$}%
   \kern-.025em\copy0\kern-\wd0
   \kern.05em\copy0\kern-\wd0
   \kern-.025em\raise.0433em\box0 }

\def\eop{{\hfill $\blacksquare$}}

\def\dt     {{\Delta t}}
\def\dN     {{\Delta N}}

\def\e{{\epsilon}}

\def\eop{{\hfill $\blacksquare$}}

\def\la {{{\lambda}}}

\def\e {{\epsilon}}

\begin {document}
\title{Nonstandard second-order formulation of the LWR model} 
\author{Wen-Long Jin\footnote{Department of Civil and Environmental Engineering, California Institute for Telecommunications and Information Technology, Institute of Transportation Studies, 4000 Anteater Instruction and Research Bldg, University of California, Irvine, CA 92697-3600. Tel: 949-824-1672. Fax: 949-824-8385. Email: wjin@uci.edu. Corresponding author}}

\maketitle
\begin{abstract}
	The seminal LWR model \citep{lighthill1955lwr,richards1956lwr} has many equivalent first-order formulations in both Eulerian and Lagrangian coordinates. In this study, we present a  second-order formulation of the LWR model based on  Phillips' model \citep{Phillips1979traffic}; but the model is nonstandard with a hyperreal infinitesimal relaxation time. 	
	Since the original Phillips model is unstable with three different definitions of stability in both Eulerian and Lagrangian coordinates, we cannot use traditional methods to prove the equivalence between the second-order model, which can be considered the zero-relaxation limit of Phillips' model,  and the LWR model,  which is the equilibrium counterpart of Phillips' model. Instead, we resort to a nonstandard method based on the equivalence relationship between second-order continuum  and car-following models established in \citep{jin2016equivalence} and prove that the nonstandard model and the LWR model are equivalent, since they have the same anisotropic car-following model and stability property. We further derive conditions for the nonstandard model to be forward-traveling and collision-free, prove that the collision-free condition is consistent with but more general than the CFL condition \citep{courant1928CFL}, and demonstrate that  only anisotropic and symplectic Euler discretization methods lead to physically meaningful solutions. We numerically solve the lead-vehicle problem and show that the nonstandard second-order model has the same shock and rarefaction wave solutions as the LWR model for both Greenshields and triangular fundamental diagrams; for a non-concave fundamental diagram we show that the collision-free condition, but not the CFL condition, yields physically meaningful results. Finally we present a  correction method to eliminate negative speeds and collisions in  general second-order models, and verify the method with a numerical example. 
	Together with \citep{jin2016equivalence}, this study presents a new approach to address the two critiques on second-order continuum models in \citep{daganzo1995requiem} and can help to guide the development and discretization of more physically meaningful second-order continuum and car-following models.

\end{abstract}

{\em Keywords}: The LWR model; Phillips' model; hyperreal infinitesimal number; anisotropic and symplectic Euler discretization methods; forward-traveling and collision-free; lead-vehicle problem; correction of second-order models.

\commentout{
{\em Highlights}:
\ben
\item We present a nonstandard second-order Phillips model.
\item The original Phillips model is unstable with three stability definitions.
\item The model and the LWR model are equivalent to the same anisotropic car-following model.
\item The collision-free condition is consistent but more general than the CFL condition.
\item Only anisotropic and symplectic Euler discretization methods are physically meaningful.
\item We present a nonstandard method to correct general second-order models.
\item This study addresses Daganzo's two critiques on second-order continuum models.
\een
}

\section{Introduction}

The seminal LWR model describes the evolution of traffic density, $k(t,x)$, speed, $v(t,x)$, and flow-rate, $q(t,x)$ in a time-space domain \citep{lighthill1955lwr,richards1956lwr} and can be derived from the following rules (Hereafter we omit $(t,x)$):
\bi
\item[R1.] The fundamental law of continuum media: $q=k v$.
\item[R2.] The fundamental diagram \citep{greenshields1935capacity}: a speed-density relation $v= V(k)$, and a flow-density relation 
\bqn
q=\phi (k)\equiv k \eta(k). \label{fd}
\eqn
Here R2 is determined by the static characteristics of traffic flow. Various fundamental diagrams have been proposed. For examples, in the triangular fundamental diagram \citep{munjal1971multilane,haberman1977model,newell1993sim}, we have:
\bqn
\eta(k)=\min\{V, W(\frac{K}k-1)\}, \label{tri-fd}
\eqn
where $V$ is the free-flow speed, $-W$ characteristic wave speed in congested traffic,  $K$ the jam density, and $\tau\equiv \frac1 {KW}$ is the time gap; in the Greenshields fundamental diagram \citep{greenshields1935capacity}, we have 
\bqn
\eta(k)&=&V (1-\frac k {K}). \label{greenshields-fd}
\eqn
\item[R3.] Conservation of vehicles: 
\bqn
\pd{k} t+\pd{kv}x=0. \label{conser_eqn}
\eqn
\ei
From R1-R3, the LWR model can be written as
\bqn
\pd{k}t+\pd{k \eta(k)}x&=&0, \label{E-S}
\eqn 
which is a scalar hyperbolic conservation law and well-defined with the following two additional implicit rules.
\bi
\item[R4.] Weak solutions: discontinuous shock waves can develop from smooth initial conditions.
\item[R5.] Entropy conditions: certain physical laws are used to pick out unique weak solutions. Traditionally the Lax or Oleinik entropy condition is used to pick out unique solutions \citep{lighthill1955lwr,ansorge1990entropy}.
\ei
Numerically the Godunov method has been successfully applied to solve the LWR model and extended as the Cell Transmission Model for network traffic \citep{daganzo1995ctm,lebacque1996godunov}.
Recently, equivalent formulations of the LWR model have been derived with different state variables and coordinate systems. For example, with the cumulative flow as the state variable in Eulerian coordinates, the LWR model is equivalent to a Hamilton-Jacobi equation, which can be solved with the minimum principle in \citep{newell1993sim} or the variational principle in \citep{daganzo2005variationalKW}. In \citep{daganzo2006ca}, the LWR model was shown to be equivalent to various car-following and cellular automaton models. In \citep{leclercq2007lagrangian}, both Hamilton-Jacobi and hyperbolic conservation law formulations of the LWR model were derived and solved in Lagrangian coordinates. In \citep{laval2013hamilton}, three equivalent Hamilton-Jacobi equations were derived and solved with the Hopf-Lax formula.

\commentout{
The LWR model and its various formulations have been successfully applied to solve various traffic problems, including the Riemann problem  \citep{lax1972shock,lebacque1996godunov}, general initial-boundary value problems \citep{newell1993sim,lebacque2003intersection,daganzo2005variationalKW}, moving bottleneck problem \citep{newell1998bottleneck}, multi-class traffic problem \citep{wong2002multiclass}, and lane-changing problem \citep{jin2010lc,jin2013multi}. It has also be extended for studying traffic dynamics in road networks through the Cell Transmission Model \citep{daganzo1995ctm} and the Link Transmission Model \citep{yperman2006mcl}.
}

It is well known that the LWR model cannot capture  scatter and hysteresis  in speed- and flow-density relations,  spontaneous stop-and-go traffic caused by instability, or bounded acceleration. These limitations have motivated many extensions of the LWR model, including second-order models  \citep[e.g.][]{payne1971PW,whitham1974PW,zhang1998theory,aw2000arz,zhang2002arz,lebacque2003acceleration}. In these models,  the conservation equation, \refe{conser_eqn}, is complemented by an additional equation for the acceleration process.  For examples, 
the acceleration equation in Zhang's model \citep{delcastillo1994reaction,zhang1998theory} can be written as
\bqn
v_t+vv_x&=&\frac{\eta(k)-v} T - (\eta'(k))^2 k k_x, \label{zhangmodel}
\eqn
where $T$ is the relaxation time. Such a second-order model reduces to the LWR model, \refe{E-S}, in equilibrium (or steady) states and is traditionally studied as a system of hyperbolic conservation laws with relaxation. It was shown in \citep{li2000global} that Zhang's model is stable and converges to the LWR model, when the relaxation time $T$ converges to zero \citep{liu1987relaxation}. That is, the LWR model is the zero-relaxation limit of Zhang's model.

In this study, we present a nonstandard second-order formulation of the LWR model based on Phillips' model \citep{Phillips1979traffic}, in which the acceleration equation without the anticipation term is simpler than Zhang's model:
\bqn
v_t+vv_x&=&\frac{\eta(k)-v} T. \label{phillipsmodel}
\eqn
In particular, we replace the relaxation time by an infinitesimal number, which equals the time step-size in numerical solutions. The formulation is nonstandard, since the infinitesimal relaxation time is  a hyperreal number in the nonstandard analysis \citep{robinson1996nonstandard}. The second-order formulation offers a number of advantages over the LWR model: non-equilibrium initial data can be incorporated, and the acceleration rate can be explicitly calculated, even though it can be a hyperreal number. The second advantage is especially helpful if one introduces bounded acceleration rates into the model. 

Furthermore, the method that we employ to analyze and solve the nonstandard second-order model is also nonstandard. In particular, we resort to the equivalence relationship between second-order continuum and car-following models established in \citep{jin2016equivalence}: we first convert a second-order continuum model into a time- and vehicle-continuous car-following model in Lagrangian coordinates, and then apply the anisotropic method to discretize derivatives in vehicles, and a symplectic Euler method to discretize the acceleration rate and speed in time. Then with their equivalent car-following models, we establish the equivalence between the nonstandard second-order model and the LWR model. We also analytically and numerically solve various lead-vehicle problems, in which the leading vehicle travels at a constant speed, to examine other properties of the model; in particular, we are interested in whether the model is stable, forward-traveling, or collision-free.

The choice of the nonstandard method is prompted by the failure of standard methods for analyzing and solving second-order models. 
\bi
\item
First, even though the nonstandard second-order formulation can be considered a zero-relaxation limit of Phillips' model, the traditional approach to proving the equivalence between the zero-relaxation limit (the nonstandard second-order formulation) and the equilibrium model (the LWR model) fails to apply, since Phillips' model is always unstable, as will be proved in Section 2.2. In contrast, with the nonstandard method, we can demonstrate that the nonstandard second-order model and the LWR model have the same car-following model and, therefore, are equivalent.  
\item
Second, traditionally second-order continuum models have been analyzed and solved as systems of hyperbolic conservation laws with relaxation \citep{liu1987relaxation}. But this approach has led to much confusion regarding the validity of second-order continuum models. For example, in many second-order continuum models, the characteristic wave speeds can be higher than vehicles' speeds; this has led to the conclusion that such models are not anisotropic, as information can travel faster than vehicles \citep{daganzo1995requiem}.  However, it was shown in \citep{leveque2001night} that, in the LWR model for night traffic,  the characteristic wave speed can also be larger than vehicles' speeds, but anisotropic solutions still exist if one converts the LWR model into a car-following model, except that such anisotropic solutions are different from those obtained with the traditional Oleinik entropy condition for hyperbolic conservation laws. Therefore, we cannot simply conclude that the LWR model with a non-concave fundamental diagram is not anisotropic; rather, 
when the fundamental diagram is non-concave, we cannot apply the traditional methods to analyze and solve the LWR model as a hyperbolic conservation law. That is, whether  characteristic wave speeds are larger than vehicles' speeds is irrelevant to the anisotropy of the LWR model, and the traditional Oleinik and Lax entropy conditions cannot be used to pick out unique, physical solutions. Since higher-order models have non-equilibrium solutions and much more complicated flow-density relations, we expect that the traditional methods are even less applicable. In contrast, with the nonstandard method,  we convert the nonstandard second-order model into a car-following model, as in \citep{leveque2001night}, and solve the latter to obtain anisotropic solutions for the former. \footnote{Many efforts have been devoted to addressing Daganzo's first critique regarding anisotropy by developing new models with the characteristic wave speeds not greater than vehicles' speeds \citep[e.g.][]{aw2000arz,zhang2002arz,jiang2002JWZ}. But we can see that the critique and the ensuing efforts to address it are largely misplaced, since we can obtain anisotropic solutions for all second-order continuum models through their equivalent car-following models, regardless whether their characteristic wave speeds are greater than vehicles' speeds or not.}

\item Third, in \citep{daganzo1995requiem}, another issue was raised for second-order models: In an extreme lead-vehicle problem when the leading vehicle is stopped, the following vehicles could travel backwards, and solutions of speeds could be negative. To the best of our knowledge,  there is no systematic method to address this critique. In contrast,  with the nonstandard method, we will be able to  derive conditions to guarantee forward-traveling and collision-free solutions in the car-following model for a second-order continuum model.
\ei

The remainder of the paper is organized as follows. In Section 2, we present the nonstandard second-order model and prove that Phillips' model is unstable. In Section 3, we show that the nonstandard second-order model and the LWR model have the same car-following formulation and stability property. In Section 4, we discuss conditions for the nonstandard second-order model to be forward-traveling and collision-free and demonstrate that the anisotropic and symplectic Euler methods are the only physically meaningful discretization methods. In Section 5, we numerically solve the lead-vehicle problem with various fundamental diagrams to verify analytical results. In Section 6, we present a method to correct existing second-order models to eliminate negative speeds and collisions. In Section 7, we conclude the study with further discussions.

\section{A nonstandard second-order model}
In this section, we first present a nonstandard Phillips' model and then demonstrate that Phillips' model is unstable, both in Eulerian and Lagrangian coordinates. The instability of Phillips' model rules out the traditional method for proving the equivalence between the nonstandard model and the LWR model.

\subsection{Nonstandard Phillips' model}
We replace the relaxation time in Phillips' model by an infinitesimal number, $\epsilon$, and obtain the following model:
\bqn
v_t+vv_x&=& \frac{\eta(k)-v}{\e}, \label{ns-lwr}
\eqn
where $\frac 1 \e$ is an infinite number. Thus the acceleration rates can be infinite in this model.

Even though infinitesimal and infinite numbers have been used ever since Archimedes' time,  especially since Leibniz used infinitesimal numbers in calculus, a rigorous treatment of such numbers was only established after the introduction of hyperreal numbers and nonstandard analysis in 1960's by \citet{robinson1996nonstandard}. The set of hyperreal numbers is an extension of the set of real numbers by including infinite and infinitesimal numbers.
That is, there exists a hyperreal (infinite) number $\Omega$, which is larger than any real number; correspondingly,  $\epsilon=\frac 1\Omega$ is a hyperreal (infinitesimal) number and smaller than any positive real number, but greater than $0$; and hyperreal numbers can be ordered and satisfy the same addition and multiplication rules for real numbers. For examples, $\Omega$ and $\epsilon$ can be defined as sequences \citep{goldblatt1998lectures}: $\Omega=\langle 1,2,3,\cdots \rangle =\langle n:n\in \mathbb{N} \rangle$, and $\epsilon=\langle 1,\frac 12,\frac 13,\cdots \rangle =\langle \frac 1n:n\in \mathbb{N} \rangle$ where $\mathbb{N}$ is the set of natural numbers.  Hyperreal numbers still form an ordered field as real numbers.

It has been shown that all mathematical arguments in nonstandard analysis can be established in standard analysis.
Thus nonstandard analysis is relatively new and not widely studied in the mathematics area. However, infinitesimal/infinite numbers and nonstandard analysis can be very useful for modeling purposes and have been widely used in physics and other areas. 
For examples, hyperreal numbers have led to simple representations of discontinuous functions, including the Dirac delta function. More importantly they enable the development of nonstandard analysis through intuitive and rigorous interpretations of derivatives and integrals with infinitesimals, which were initiated by Leibniz \citep{davis2005applied}.
In \citep{hanqiao1986applications,vandenberg1998relation}, nonstandard analysis was applied to solve the heat equation under initial conditions given by a Dirac delta function. In these studies, $\epsilon$ is defined as the limit of the time-step size, $\dt$, and it was shown that $\epsilon$ can be simply replaced by $\dt$ in the discrete version of a differential equation and the continuous and discrete equations are equivalent for an infinitesimal $\dt$.

In this study, we also choose $\e$ such that it equals the time-step size, $\dt$, in the discrete version. In a sense, $\e$ is a limit:
\bqn
\e&=&\lim_{\dt\to 0^+} \dt,
\eqn
which is an indefinitely small number smaller than any fixed small positive number.  In \citep{jin2015ltm}, an indicator function $H(y)$ for $y\geq 0$ was defined as
\bqs
H(y)&=&\lim_{\dt\to 0^+} \frac y \dt=\cas{{ll}0, & y=0\\+\infty, &y>0}
\eqs
and used in defining both demand and supply functions for the link transmission model. Such an indicator function was also used in formulating the point queue models \citep{jin2015pqm} and defining a Dirac delta function in \citep{jin2016instantaneous}. Thus, $H(y)$ can be represented by the hyperreal number $\epsilon$ as
\bqs
H(y)&=&\frac y\epsilon,
\eqs
or  \refe{ns-lwr} can be re-written with the indicator function
\bqs
v_t+vv_x&=& H(\eta(k)-v).
\eqs

\subsection{Instability of Phillips' model}
Phillips' model, \refe{conser_eqn} and \refe{phillipsmodel}, is a system of hyperbolic conservation laws with relaxation, and the LWR model is its equilibrium counterpart. For the system of hyperbolic conservation laws, the two characteristic wave speeds are always identical: $\la_1(k,v)=\la_2(k,v)=v$. Therefore, it is non-strictly hyperbolic. We denote the characteristic wave speed of the LWR model by $\la_*(k)=\phi'(k)=\eta(k)+k\eta'(k)$.

Since nonstandard Phillips' model can be considered as the zero-relaxation limit of Phillips' model, if Phillips' model is stable, then the equivalence between its zero-relaxation limit and the LWR model can be established by following the traditional methods \citep{liu1987relaxation,li2000global}. But in this subsection we demonstrate that Phillips' model is unstable with three different definitions of stability in both Eulerian and Lagrangian coordinates.

\begin{theorem} For Phillips' model, its approximate viscous LWR model can be written as
\bqn
k_t +(k \eta(k))_x&\approx&- T \pd{}x \left( (\la_*(k)-\eta(k))^2 k_x\right), \label{diffusioneqn}
\eqn	
which is diffusively unstable near an equilibrium state.
\end{theorem}
{\em Proof}. Following \citep{liu1987relaxation}, we apply the Chapman-Enskog expansion for a small perturbation around an equilibrium state: $v=\eta(k)+v_1$. Here we assume that both $v_1$ and its derivatives are small; i.e., $v_1>>\pd{v_1}x >> \frac{ \partial ^2 v_1}{\partial x^2}$ holds for diffusive waves as $t\to \infty$.

Thus \refe{phillipsmodel} can be re-written as
\bqs
\pd{}t(\eta(k)+v_1) +\pd{} x \left(\frac 12 (\eta(k)+v_1)^2\right) = -\frac{v_1}T,
\eqs
which can be approximated by
\bqs
v_1&\approx&-T \left[\pd{}t \eta(k) +\pd{} x \frac 12 (\eta(k))^2\right] =-T \eta'(k) (k_t+\eta(k) k_x),
\eqs
by omitting the derivatives of $v_1$.

Along the characteristic wave direction, we have $\pd{}t+\phi'(k) \pd{}x \approx 0$, which leads to $k_t\approx -(\eta(k)+k \eta'(k)) k_x$. Thus we have
\bqs
v_1&\approx&  T [\eta'(k)]^2 kk_x,
\eqs

Substituting this into \refe{conser_eqn}, we have
\bqs
k_t+(k\eta(k))_x&\approx &- T \pd {}x ((k \eta'(k))^2 k_{x}),
\eqs
which is equivalent to \refe{diffusioneqn}. Since the coefficient of the viscous term is negative, the diffusion process is unstable. 
\eop

\begin{theorem} Phillips' model is linearly unstable. 	
	\end{theorem}
	{\em Proof}.
We follow \citep[][Section 3.1]{whitham1974PW} to prove the linear instability of Phillips' model.
We first linearize Phillips' model  about an equilibrium state at $k_0$ and $v_0$, which satisfy $v_0=\eta(k_0)$, by assuming $k=k_0+\kappa$ and $v=v_0+\nu$, where  $\kappa(t,x)$ and $\nu(t,x)$ are small perturbations. Omitting higher-order terms, we have 
\bqs
\kappa_t +v_0 \kappa_x +k_0 \nu_x &=&0,\\
\nu_t +v_0 \nu_x &=&\frac{\eta'(k_0)\kappa-\nu}T.
\eqs 
Then we take partial derivatives of the second equation with respect to $x$ and substitute $\nu_x$ from the first equation to obtain
\bqs
T \left(\pd{}{t}+v_0 \pd{}{x}\right)^2 \kappa&=&-k_0 \eta'(k_0) \pd \kappa x- \left(\pd{}{t}+v_0 \pd{}{x}\right) \kappa.
\eqs

We consider the exponential solution $\kappa=e^{i(m x- \omega t)}$, where $m$ is a real number. The exponential solutions are stable if and only if the imaginary part of $\omega$ is negative. 
Then the above linearized system is equivalent to
\bqn
-T(\omega+mv_0)^2=-i k_0\eta'(k_0) m +i (\omega+mv_0).
\eqn
The above equation can be rewritten as  $\omega^2+(2b_1+i 2b_2) \omega +(d_1+i d_2)=0$ with
\bqs
2b_1&=&2mv_0, \quad 2b_2=\frac 1T,\\
d_1&=&m^2 v_0^2, \quad d_2=\frac{v_0-k_0\eta'(k_0)}T m.
\eqs
According to \citep{abeyaratne2014macroscopic}, both roots of the equation have negative imaginary parts if and only if $b_2>0$ and $4b_1b_2 d_2-4d_1 b_2^2>d_2^2$, or equivalently
\bqs
T&>&0,\\
(k_0\eta'(k_0))^2&<&0,
\eqs
which is impossible. Therefore Phillips' model is linearly unstable.
\eop 

 As shown in \citep{makigami1971traffic}, the evolution of a traffic stream can be captured by a surface in a three-dimensional space of time $t$, location $x$, and vehicle number $N$. Here the positive direction of $x$ is the same as the traffic direction, and $N$ is the cumulative number of vehicles passing location $x$ at $t$ after a reference vehicle \citep{moskowitz1965discussion}. Therefore, $N$ increases in $t$ but decreases in $x$; i.e., a leader's number is smaller than the follower's.
 Among these primary three variables, two of them are independent and the other dependent. That is, there exist the following functions:
 \bqn
 N&=&n(t,x),\\
 x&=&X(t,N).
 \eqn
Here $(t,x)$ forms the Eulerian coordinates, and $(t,N)$ the Lagrangian coordinates.

In Eulerian coordinates, $k=-n_x$, and $v=-n_t/n_x$. Transforming all the variables from the Eulerian coordinates into the Lagrangian coordinates, we then have $v=X_t$, $v_t+vv_x=X_{tt}$, and $k=-\frac 1{X_N}$. Thus Phillips' model can be rewritten as a time- and vehicle-continuous car-following model in the Lagrangian coordinates as
\bqn
X_{tt}&=&\frac{\theta(-X_N)-X_t} T, \label{L-P-phillips}
\eqn
where 
\bqn
\theta(S)=\eta\left(\frac 1 S\right)
\eqn
 is the speed-spacing relation. For \refe{L-P-phillips}, the time-continuous and vehicle-discrete car-following model can be written as
\bqn
X_{tt}(t,N)&=&\frac 1T \left( \theta\left(\frac{X(t,N-\dN)-X(t,N)}\dN\right)-X_t(t,N)\right), \label{L-P-vdis}
\eqn
which is the optimal velocity model with an arbitrary $\dN$ \citep{bando1995dmt}. Note that, since traffic flow is usually anisotropic and information propagates from the leading vehicle $N-\dN$ to the following $\dN$, we approximate $-X_N(t,N)$ by $\frac{X(t,N-\dN)-X(t,N)}\dN$ \citep{leveque2001night,jin2016equivalence}.

\begin{theorem} \label{thm:phillipslinear}
	Phillips' model in the Lagrangian coordinates, \refe{L-P-phillips} or \refe{L-P-vdis}, is string unstable. That is,  a small disturbance in a leading vehicle's speed gets amplified along a traffic stream.
	\end{theorem}
{\em Proof}.
We assume that there are small monochromatic disturbances to vehicles' speeds near an equilibrium state with a spacing of $s_0$ and a speed of $v_0$, where $s_0=\theta(s_0)$. Here we denote
$X_t(t,N-j \dN)=v_0+\nu_{N-j \dN} e^{i\omega t}$ for $j=0,1,\cdots, \frac 1\dN$. 
Then $X_{tt}(t,N)=i\omega \nu_N e^{i\omega t}$, and  $X(t,N-j\dN)=v_0t+\frac{\nu_{N-j\dN}}{i\omega} e^{i\omega t} +X(0,N-j\dN)$. Thus, $X(t,N-\dN)-X(t,N)=s_0\dN+\frac{\nu_{N-\dN}-\nu_N}{i\omega} e^{i\omega t}$, where $s_0=X(0,N-1)-X(0,N)$ is the initial equilibrium spacing. From \refe{L-P-vdis}, we have
\bqs
i\omega \nu_N e^{i\omega t} &=& \frac 1T \left( \theta\left(s_0+\frac{\nu_{N-\dN}-\nu_N}{i\omega \dN} e^{i\omega t}\right)- v_0-\nu_N e^{i\omega t} \right).
\eqs
With the Taylor series expansion of the right-hand side, we obtain
\bqs
i\omega \nu_N e^{i\omega t} &\approx&  \frac 1T\left(\theta'(s_0) \frac{\nu_{N-\dN}-\nu_N}{i\omega\dN} e^{i\omega t}- \nu_N e^{i\omega t}\right),
\eqs
which leads to
\bqs
i\omega \nu_N  &\approx& \frac1T \left(\theta'(s_0) \frac{\nu_{N-\dN}-\nu_N}{i\omega\dN}- \nu_N \right).
\eqs
Thus we have for small $\dN$
\bqs
\frac{\nu_N} {\nu_{N-\dN}}&\approx& \frac{\theta'(s_0)}{\theta'(s_0)-T\omega^2 \dN+i\omega \dN} \approx 1-\frac 1{\theta'(s_0)} (i\omega-T\omega^2) \dN,
\eqs
which leads to by omitting the terms with $\dN^2$
\bqs
\left|\frac{\nu_N} {\nu_{N-\dN}}\right|^2\approx 1+ \frac{2T\omega^2 \dN }{\theta'(s_0)}.
\eqs
Thus we have
\bqs 
\left|\frac{\nu_N} {\nu_{N-1}}\right| = \lim_{\dN\to 0} \left|\frac{\nu_N} {\nu_{N-\dN}}\right|^{\frac 1\dN} \approx e^{ \frac{T\omega^2  }{\theta'(s_0)}}.
\eqs
Since generally speed increases in spacing; i.e., $\theta'(s_0)>0$, $|\nu_N|>|\nu_{N-1}|$. That is, the following vehicle's speed oscillation magnitude is always larger than the leader's. Therefore a small disturbance in a leading vehicle's speed gets amplified along a traffic stream, and Phillips' model in the Lagrangian coordinates, \refe{L-P-phillips}, is string unstable. 
\eop

\section{Equivalence and stability}
In steady states when the acceleration rate is zero, $v=\eta(k)$ for \refe{conser_eqn} and \refe{ns-lwr}, and the nonstandard second-order model is equivalent to the LWR model. 
In this section we prove the equivalence between the nonstandard second-order model and the LWR model under general conditions. We also discuss its stability property.

\subsection{Equivalence between the nonstandard second-order model and the LWR model}
Here we apply the conversion method developed in \citep{jin2016equivalence}: we first convert the nonstandard model into equivalent time- and vehicle-continuous car-following model in Lagrangian coordinates, and then discretize the latter into time- and/or vehicle-discrete car-following models. 

First, the nonstandard second-order model,  \refe{conser_eqn} and \refe{ns-lwr}, is equivalent to the following time- and vehicle-continuous car-following model in Lagrangian coordinates:
\bqn
X_{tt}&=&\frac{\theta(-X_N)-X_t} \e.  \label{con-cf}
\eqn
Then we discretize $X_N$ in a backward fashion:
\bqn
X_N(t,N)&=& \frac{X(t,N)-X(t,N-\dN)}\dN. \label{backward-difference}
\eqn
We refer to \refe{backward-difference} as an anisotropic method, since, with such a discretization in vehicles, the above model is automatically anisotropic; i.e., the acceleration rate of vehicle $N$ is only impacted by its leading vehicle $N-\dN$, but not its following vehicle $N+\dN$.
Thus the time-continuous and vehicle-discrete car-following model of \refe{con-cf} is 
\bqn
X_{tt}(t,N)&=&\frac 1\e \left(\theta\left(\frac{X(t,N-\dN)-X(t,N)}\dN\right)-X_t(t,N)\right). \label{ns-cf}
\eqn

Furthermore, for a discrete time axis with a time step-size of $\dt$, both $X_{tt}(t,N)$ and $X(t,N)$ are defined at $t$, but the speed $X_t(t,N)$ is actually defined at $t-\frac 12 \dt$, since it approximates the average speed between $t$ and $t-\dt$. Therefore we apply the following explicit and implicit Euler discretization methods for the acceleration rate and speed, respectively, as follows \citep{jin2016equivalence}:
\bsq\label{explicit-implicit}
\bqn
X_{tt}(t,N)&=&\frac{X_t(t+\dt,N)-X_t(t,N)}{\dt},\\
X_t(t+\dt,N)&=&\frac{X(t+\dt,N)-X(t,N)}\dt.
\eqn
\esq
Note that \refe{explicit-implicit} is a symplectic Euler method, if we re-write the time-continuous and vehicle-discrete car-following model, \refe{ns-cf}, as two Hamiltonian equations for $X(t,N)$ and $X_t(t,N)$ \citep[][Chapter 6]{hairer2006geometric}.
Setting $\epsilon=\dt$, we then obtain the time- and vehicle-discrete car-following model from \refe{ns-cf}:
\bsq
\bqn
X_t(t+\dt,N)&=& \theta\left(\frac{X(t,N-\Delta N)-X(t,N)}{\Delta N} \right), \label{ns-dis1} \\
X(t+\dt,N)&=&X(t,N)+\dt \cdot \theta\left(\frac{X(t,N-\Delta N)-X(t,N)}{\Delta N} \right).\label{ns-dis2}
\eqn
\esq

We can see that \refe{ns-dis2} is the same as Equation (12) in \citep{leclercq2007lagrangian}, derived from the Godunov scheme of the LWR model. Therefore, the nonstandard second-order model, \refe{conser_eqn} and \refe{ns-lwr},  is equivalent to the LWR model, \refe{E-S}.
It is essentially a first-order model but with hyperreal acceleration rates.

\subsection{Stability}

\begin{theorem}
	The time-continuous and vehicle-discrete car-following model, \refe{ns-cf}, is linearly stable. That is, a small disturbance in a leading vehicle's speed is not amplified along a traffic stream.
\end{theorem}
{\em Proof}. Following the proof of Theorem \ref{thm:phillipslinear}, we have  for small $\dN$ and $\e$
\bqs
\frac{\nu_N} {\nu_{N-\dN}}&=& \frac{\theta'(s_0)}{\theta'(s_0)-\e\omega^2 \dN+i\omega \dN} \approx 1-\frac 1{\theta'(s_0)} i\omega \dN,
\eqs
since $\e$ is the same order of $\dN$. 
Thus we have
\bqs
\left|\frac{\nu_N} {\nu_{N-\dN}}\right|^2\approx 1 -\frac{\omega^2 \dN^2}{(\theta'(s_0))^2}<1,
\eqs
and
\bqs 
\left|\frac{\nu_N} {\nu_{N-1}}\right| = \lim_{\dN\to 0} \left|\frac{\nu_N} {\nu_{N-\dN}}\right|^{\frac 1\dN}= 1.
\eqs

Therefore a small disturbance in a leading vehicle's speed is not amplified along a traffic stream, and the nonstandard car-following model, \refe{ns-cf}, is stable. 
\eop

The stability property of the car-following model is the same as that of the original LWR model.

\section{Forward-traveling and collision-free properties}
In the preceding section, we demonstrated that, for the nonstandard second-order model, \refe{conser_eqn} and \refe{ns-lwr}, there exists a unique speed-density relation, its equivalent car-following models, \refe{ns-cf}, \refe{ns-dis1}, and \refe{ns-dis2}, admit anisotropic solutions, and it is  stable. 

In this section, we discuss two additional properties of the model: forward-traveling and collision-free. We call a model forward-traveling, if vehicles' speeds are always non-negative. We call a model collision-free, if vehicles' spacings are not smaller than the jam spacing,  $S\equiv \frac 1K$; i.e., traffic densities  are not greater than the jam density. The forward-traveling property was initially introduce in \citep{daganzo1995requiem} as a criterion to determine the validity of a second-order continuum model. In this study we extend this criterion for car-following models.

\subsection{Conditions for forward-traveling and collision-free}
In this subsection, we derive conditions for the time- and vehicle-discrete car-following models, \refe{ns-dis1} and \refe{ns-dis2}, to be forward-traveling and collision-free.

\begin{theorem} \label{thm:nonnegative} The speed in \refe{ns-dis1} is  non-negative, if the spacing is not smaller than the jam spacing; i.e., if $X(t,N-\Delta N)-X(t,N) \geq S \Delta N  $. That is, both  \refe{ns-dis1} and \refe{ns-dis2} are forward-traveling if they are collision-free.
	\end{theorem}
{\em Proof}. For a valid speed-spacing relation, $v=\theta(s)$, the speed is non-negative, if the spacing is not smaller than the jam spacing. Therefore, if $X(t,N-\Delta N)-X(t,N) \geq S \Delta N$, the speed is non-negative. That is, both  \refe{ns-dis1} and \refe{ns-dis2} are forward-traveling if they are collision-free. \eop

Note that, for other car-following models, forward-traveling and collision-free properties may be independent of each other. In the following theorem, we derive the sufficient and necessary condition for \refe{ns-dis1} and \refe{ns-dis2} to be collision-free.

\begin{theorem} \refe{ns-dis2}  is collision-free, if and only if  
$\dt$ and $\dN$ satisfy the following condition:
\bqn
\frac{\dN}\dt &\geq &\max_{k\in[0,K]} \frac{\phi(k)}{1-\frac kK}. \label{collision-free-cond}
\eqn
\end{theorem}
{\em Proof.} Assuming that the model is collision-free at $t$; i.e., $X(t,N-\dN)-X(t,N)\geq S \dN$. At $t+\dt$,  we have
\bqs
X(t+\dt,N-\dN)-X(t+\dt,N)&\geq&  X(t,N-\dN)-X(t+\dt,N),
\eqs
where the equal sign holds when the leader, vehicle $N-\dN$, suddenly stops. From  \refe{ns-dis2}, we have 
\bqs
X(t,N-\dN)-X(t+\dt,N)&=&X(t,N-\dN)-X(t,N)\\
&&
-\dt \cdot \theta\left(\frac{X(t,N-\Delta N)-X(t,N)}{\Delta N} \right).
\eqs
Denoting $k=\frac {\dN}{X(t,N-\dN)-X(t,N)}$, we have
\bqs
X(t,N-\dN)-X(t+\dt,N)&=& \frac\dN k-\dt \eta(k).
\eqs
Thus the model is always collision-free, if and only if for any $k\in [0,K]$
\bqs
\frac\dN k-\dt \eta(k) \geq S \dN =\frac \dN K,
\eqs
which is equivalent to \refe{collision-free-cond}.
\eop

In \citep{leclercq2007lagrangian},  it was shown that  \refe{ns-dis2} is a Godunov finite difference equation for a hyperbolic conservation formulation of the LWR model in Lagrangian coordinates. For \refe{ns-dis2}  to be stable, consistent, and, therefore, convergent, the traditional CFL condition \citep{courant1928CFL} was derived as 
\bqn
\frac{\dN}\dt &\geq & \max_{k\in[0,K]} | -\eta'(k) k^2 |. \label{cfl-cond}
\eqn
Clearly the collision-free condition, \refe{collision-free-cond}, is not the same as the CFL condition, \refe{cfl-cond}. But we demonstrate that they are equivalent, when the speed-density relation is non-increasing, and the flow-density relation is concave.

\begin{theorem} The collision-free condition, \refe{collision-free-cond}, is equivalent to the CFL condition, \refe{cfl-cond}, when (i) the speed-density relation is non-increasing; i.e., $\eta'(k)\leq 0$;  and (ii) for $k<K$
\bqn
k\eta''(k)+2 \eta'(k) \leq 0, \label{cfl-equivalence}
\eqn
which is equivalent to the concavity condition for the flow-density relation:
\bqn
\phi''(k)&\leq &0. \label{concavity-cond}
\eqn
In particular, the two conditions are equivalent for the Greenshields and triangular fundamental diagrams.  
\end{theorem}
{\em Proof}. We first assume that the speed-density relation is non-increasing; i.e., $\eta'(k)\leq 0$. This is generally true for many traffic systems, but not for the night traffic system considered in \citep{leveque2001night}.

We then assume that \refe{cfl-equivalence} is satisfied for $k<K$. It is straightforward to show that \refe{cfl-equivalence} is equivalent to \refe{concavity-cond}. From \refe{cfl-equivalence} we have $(K-k)(k\eta''(k)+2 \eta'(k))\leq 0$. Thus $(K-k) k \eta'(k)+K\eta(k)$ is a non-increasing function in $k$, since its derivative is $(K-k)(k\eta''(k)+2 \eta'(k))$. Thus $(K-k) k \eta'(k)+K\eta(k) \geq (K-K) K \eta'(K)+K\eta(K)=0$. Further we can see that $\frac{\phi(k)}{1-\frac kK}$ is non-decreasing function, since its derivative is 
\bqs
\frac{\phi'(k)(1-\frac kK)+\phi(k)\frac 1K}{(1-\frac kK)^2}&=&
\frac{(K-k) k \eta'(k)+K\eta(k)}{K (1-\frac kK)^2} \geq 0.
\eqs
Thus \refe{collision-free-cond} is equivalent to $\frac \dN\dt \geq  \frac{\phi(k)}{1-\frac kK}|_{k=K}= -\eta'(K) K^2$.

The CFL condition, \refe{cfl-cond}, is equivalent to
\bqs
\frac{\dN}\dt &\geq & \max_{k\in[0,K]}  -\eta'(k) k^2 .
\eqs
When \refe{cfl-equivalence} is satisfied for $k<K$, then $-\eta'(k) k^2$ is non-decreasing in $k$, and  $\max_{k\in[0,K]}  -\eta'(k) k^2= -\eta'(K) K^2$. Thus \refe{cfl-cond} is equivalent to 
$\frac \dN\dt \geq -\eta'(K) K^2$.

Therefore, \refe{collision-free-cond} and \refe{cfl-cond} are equivalent when  the speed-density relation is non-increasing, and the flow-density relation is concave.

Then we can show that the Greenshields and triangular fundamental diagrams satisfy the two conditions.
\ben
\item
For the Greenshields fundamental diagram, \refe{greenshields-fd}, we have (i) $\eta'(k)=-\frac VK <0$, and (ii) for $k<K$
\bqs
k\eta''(k)+2 \eta'(k) &=&-\frac{2V}K<0.
\eqs
Thus,\refe{collision-free-cond} and \refe{cfl-cond} are equivalent to 
\bqn
\frac \dN\dt\geq VK.
\eqn
\item For the triangular fundamental diagram, \refe{tri-fd}, we have (i) $\eta'(k)=0$ for $k\leq \frac{W}{V+W}K$ and $\eta'(k)=-\frac{WK}{k^2}<0$ for $k>\frac{W}{V+W}K$, and (ii) for $k<K$
\bqs
k\eta''(k)+2 \eta'(k) &=&0.
\eqs
Thus, \refe{collision-free-cond} and \refe{cfl-cond} are equivalent to 
\bqn
\frac \dN\dt\geq WK.
\eqn
\een
\eop

Note that the CFL condition is derived from the numerical perspective with respect to convergent methods, but the collision-free condition is derived from the physical perspective with respect to safety. 
The consistence between the collision-free and CFL conditions for normal fundamental diagrams suggests that the physical perspective is meaningful. Further in Section 5, we demonstrate that the collision-free condition is more general and still applies to non-concave fundamental diagrams. 

\subsection{Other discretization methods}
For the time- and vehicle-continuous car-following model, \refe{con-cf}, we apply the anisotropic and symplectic Euler discretization methods proposed in \citep{jin2016equivalence}: \refe{backward-difference} for discretizing the spacing in vehicles, and \refe{explicit-implicit} for discretizing the acceleration rate and speed in time. In this subsection, we examine the collision-free property with other discretization methods.

First, we apply the symplectic Euler method to discretize the acceleration and speed in time to obtain the time-discrete and vehicle-continuous form of  \refe{con-cf}:
\bqn
X(t+\dt, N)&=&X(t,N)+\dt \cdot \theta(-X_N(t,N)). \label{time-dis-cf}
\eqn
But we replace the anisotropic vehicle-discretization method, \refe{backward-difference}, by one of the following non-anisotropic discretization methods:
\ben
\item Forward difference method: 
\bqn
X_N(t,N) &=& \frac{X(t,N+\dN)-X(t,N)}\dN. \label{forward-dif}
\eqn
\item Arithmetic central difference method: 
\bqn
X_N(t,N) &=& \frac{X(t,N+\dN)-X(t,N-\dN)}{2\dN}, \label{arith-cd}
\eqn
which is the arithmetic average of the backward and forward differences.
\item Harmonic central difference method, which corresponds to the arithmetic central difference method for density \citep{leveque2001night}: 
\bqn
X_N(t,N) &=& \frac 2 { \frac{\dN}{X(t,N+\dN)-X(t,N)} +\frac {\dN}{ X(t,N) -X(t,N-\dN)}  }, \label{har-cd}
\eqn
which is the harmonic average of the backward and forward differences.
\een 
All these discretization methods are non-anisotropic, since the following vehicle's location would impact the leading vehicle's decision. In the following theorem, we demonstrate that such non-anisotropic discretization methods cannot guarantee the collision-free property of the resulted time- and vehicle-discrete car-following models.

\begin{theorem}
The time- and vehicle-discrete car-following model of \refe{time-dis-cf} with any of the non-anisotropic discretization methods, \refe{forward-dif}, \refe{arith-cd}, or \refe{har-cd}, is not collision-free, if the speed-density function is greater than zero for any spacing greater than the jam spacing. In other words, it is not safe, if a vehicle determines its speed by looking backwards.
\end{theorem}
{\em Proof}. 
We consider an extreme lead-vehicle problem, in which the leading vehicle $N-\dN$ stops at  an intersection. This corresponds to a red-light scenario. We assume that at $t$ vehicle $N-\dN$ is stopped, the spacing between vehicles $N$ and $N-\dN$ is $X(t,N-\dN)-X(t,N) = S \dN$, but the spacing between vehicles $N+\dN$ and $N$ is $X(t,N)-X(t,N+\dN) > S \dN$. So there is no collision at $t$.

But at $t+\dt$ we have from \refe{time-dis-cf}
\bqs
X(t+\dt,N-\dN)-X(t+\dt,N)&=&X(t,N-\dN)-X(t,N)-\dt\cdot \theta(-X_N(t,N))
\\&=&S \dN- \dt\cdot \theta(-X_N(t,N)).
\eqs
If we use any of the the non-anisotropic discretization methods, \refe{forward-dif}, \refe{arith-cd}, or \refe{har-cd}, then $-X_N(t,N)>S \dN$, which leads to $\theta(-X_N(t,N))>0$ for a speed-density function, in which the speed is greater than zero for any spacing greater than the jam spacing. Therefore, $X(t+\dt,N-\dN)-X(t+\dt,N)<S \dN$, and vehicle $N$ collides into vehicle $N-\dN$.
Thus the model is not collision-free. 

In addition, the model is not collision-free for  other non-anisotropic discretization methods, since $-X_N(t,N)$ would contain the upstream spacing, $\frac{X(t,N+\dN)-X(t,N)}\dN$, and is greater than the jam spacing for the red-light scenario.
\eop

From the above theorem, we can see that the  discretization method in vehicles has to be anisotropic for a model to be collision-free. Even though there can be other anisotropic discretization methods, which can involve multiple leaders, such multiple-leader car-following models are not equivalent to the LWR model.

Next, we apply the anisotropic discretization method for the spacing but use the explicit  Euler method to discretize both the acceleration rate and speed in time; i.e., we replace \refe{explicit-implicit} by
\bsq\label{explicit-explicit}
\bqn
X_{tt}(t,N)&=&\frac{X_t(t+\dt,N)-X_t(t,N)}{\dt},\\
X_t(t,N)&=&\frac{X(t+\dt,N)-X(t,N)}\dt.
\eqn
\esq
 Thus the corresponding time- and vehicle-discrete car-following model of \refe{con-cf} can be written as
\bqn
X(t+\dt,N)&=&X(t,N)+ \dt \cdot \theta\left(\frac{X(t-\dt,N-\Delta N)-X(t-\dt,N)}{\Delta N} \right). \label{cf-e-e}
\eqn
In the following theorem, we demonstrate that this model is not collision-free.

\begin{theorem} The car-following model, \refe{cf-e-e}, is not collision-free.
\end{theorem}
{\em Proof}. Again we just need to present one scenario, in which collision can occur at $t+\dt$ even though there is no collision at $t$ or $t-\dt$.

Consider a red-light scenario: vehicle $N-\dN$ is stopped since $t-\dt$, $X(t,N-\dN)-X(t,N)=S\dN$, and $X(t-\dt,N-\dN)-X(t-\dt,N)>S\dN$. Then at $t+\dt$, we have from \refe{cf-e-e}:
\bqs
X(t+\dt,N-\dN)-X(t+\dt,N)&=&X(t,N-\dN)-X(t,N)\\&&- \dt \cdot \theta\left(\frac{X(t-\dt,N-\Delta N)-X(t-\dt,N)}{\Delta N} \right)\\
&=&S\dN \\&&-\dt \cdot \theta\left(\frac{X(t-\dt,N-\Delta N)-X(t-\dt,N)}{\Delta N} \right)\\&<& S\dN.
\eqs
That is, vehicle $N$ cannot stop at the safe distance and collides into vehicle $N-\dN$ at $t+\dt$.
\eop

We can show that the Runge-Kutta and many other  discretization methods for the acceleration and speed cannot guarantee the collision-free property of the resulted car-following models. 
Therefore, the anisotropic and symplectic Euler methods proposed in \citep{jin2016equivalence}, \refe{backward-difference} and \refe{explicit-implicit}, appear to be the only physically meaningful methods to discretize the time- and vehicle-continuous car-following model, \refe{con-cf}, and  \refe{ns-dis2} the only physically meaningful car-following model that is equivalent to the nonstandard second-order model, \refe{conser_eqn} and \refe{ns-lwr}. 

\section{Numerical examples}
In this section, we numerically solve the car-following models of the nonstandard second-order model,  \refe{ns-dis1} and \refe{ns-dis2}. 
We discretize a simulation time duration $[T_0,T_0+J \dt]$ into $J$ intervals and denote $t^j=T_0+j \dt$ ($j=0,\cdots,J$). We also divide a plotoon $[0,M \dN]$ into $M$ vehicles and denote $N_m=m \dN$ ($m=0,\cdots,M$). We further denote $U_m^j=X_t(t^j,N_m)$ and $Y_m^j=X(t^j,N_m)$. Then \refe{ns-dis1} and \refe{ns-dis2} can be respectively rewritten as
\bsq \label{discrete_lp2_approx}
\bqn
U^{j+1}_m &=& \theta\left(\frac{Y^j_{m-1} -Y^j_m}{\dN}\right),\\
Y^{j+1}_m &= &Y^j_m + \dt \cdot  U^{j+1}_m . 
\eqn
Thus the acceleration rate of vehicle $m$ at $j\dt$ is
\bqn
Z_m^{j}&=&\frac {U_m^{j+1}-U_m^j}\dt.
\eqn
\esq

Here we consider the lead-vehicle problem under the following initial and boundary conditions:  the boundary conditions are determined by vehicle 0's trajectory ($j\geq 0$):
\bqs
Y^{j+1}_0 &= & Y^j_0 +\dt \cdot U^{j+1}_0,
\eqs
where $Y_0^0=0$, and the leading vehicle's speed $U^{j+1}_0=v_2$ is given and constant; the initial spacings  of all following vehicles are constant at $Y^0_m=-m s_1 \dN$ for $m=1,\cdots, M$, where the initial density is $k_1=1/s_1$, and their initial speeds are constant at $U^0_m=v_1=\eta(k_1)$. Thus the lead-vehicle problem corresponds to a Riemann problem with
\bqs
k(0,x)&=&\cas{{ll} k_1, & x<0 \\ k_2, & x>0}
\eqs
where $v_2=\eta(k_2)$. As we know, the Riemann problem of the LWR model with a concave fundamental diagram is solved by either a shock or a rarefaction wave \citep{lebacque1996godunov}. In this section we will compare the numerical results of \refe{discrete_lp2_approx} with the theoretical Riemann solutions.

\subsection{Greenshields fundamental diagram}
In this subsection we solve the lead-vehicle problem with the Greenshields fundamental diagram, \refe{greenshields-fd}, where $S=7$ m, $K=1/S$, and $V=20$ m/s. We set $\dt=0.35 \dN$, which satisfies the CFL condition and the collision-free condition.

\begin{figure} \bc
$\ba{c@{\hspace{0.1in}}c}
\includegraphics[width=3.2in]{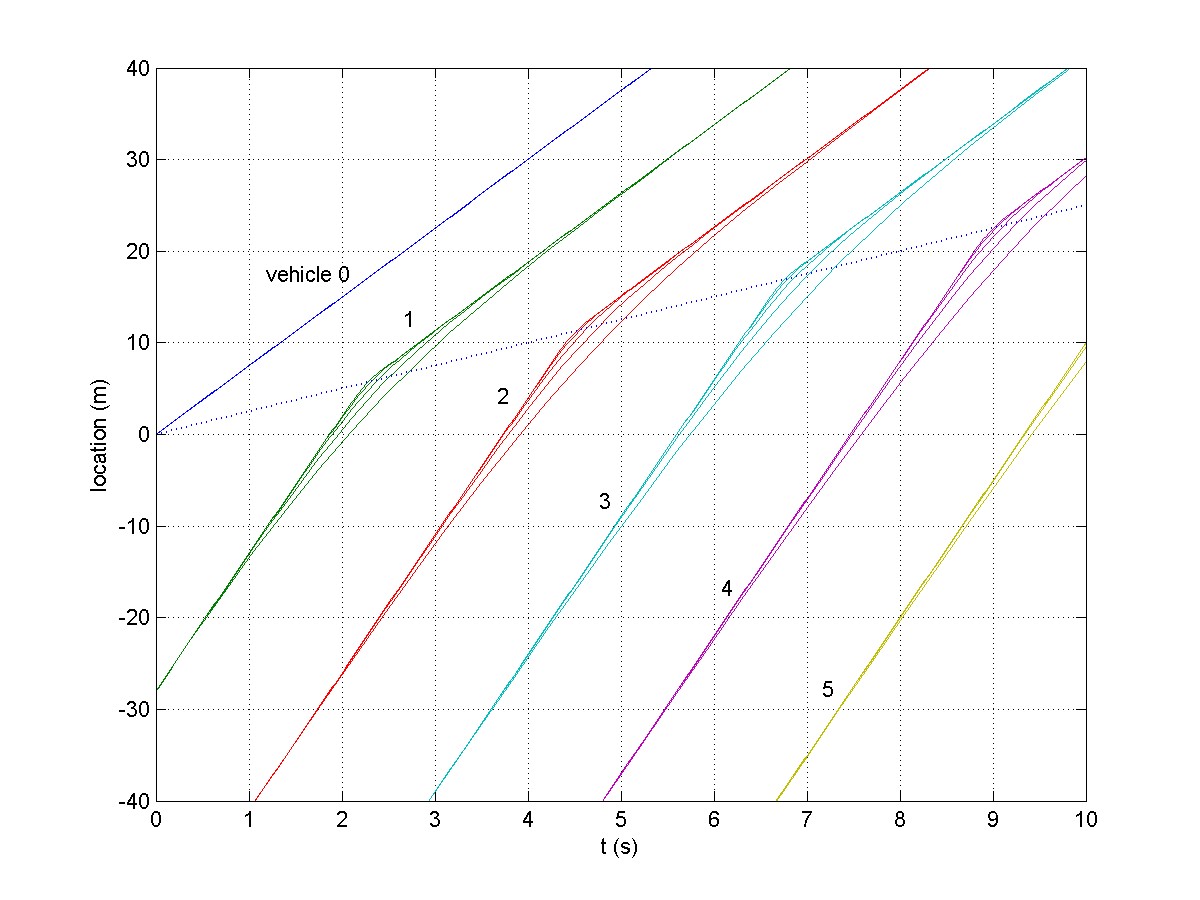} &
\includegraphics[width=3.05in]{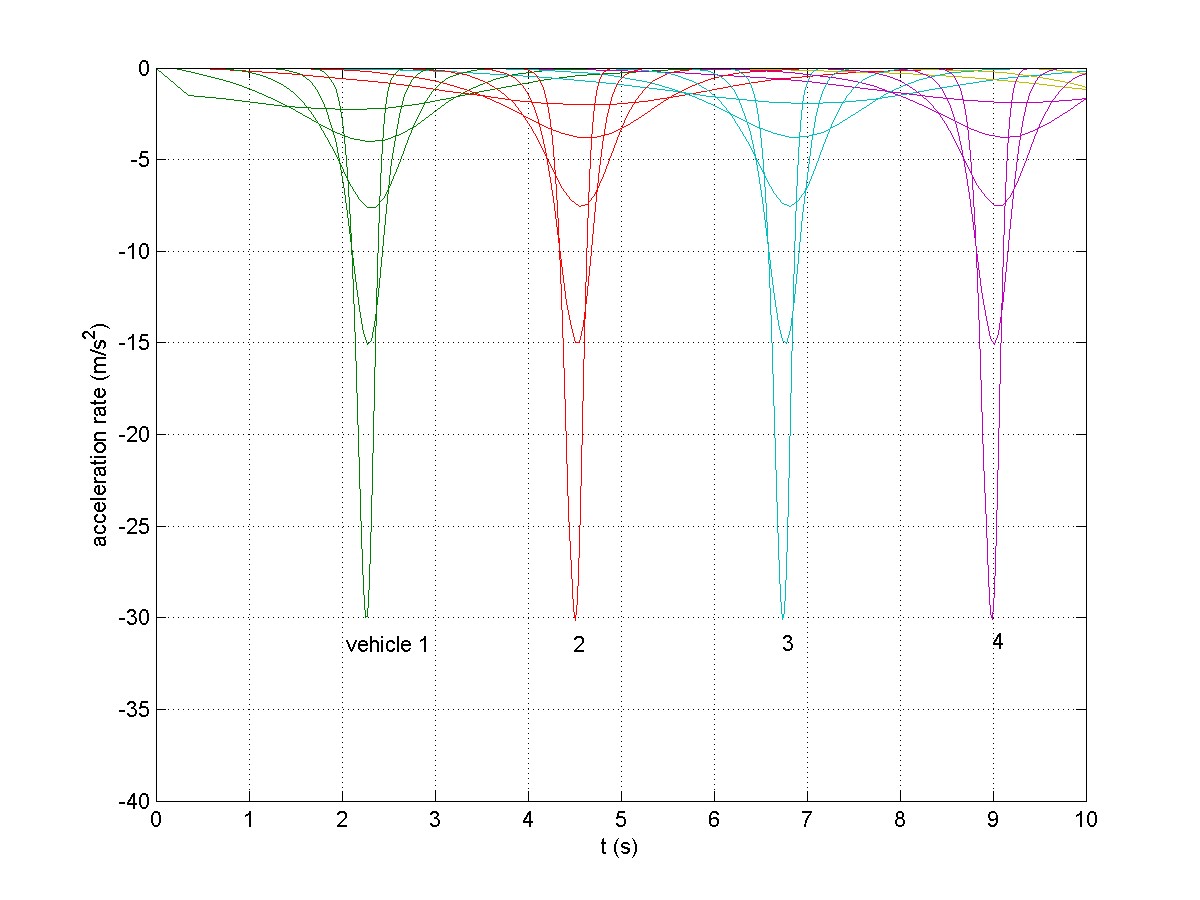} \\
\mbox{\bf (a)} &
    \mbox{\bf (b)}\\
\includegraphics[width=3.2in]{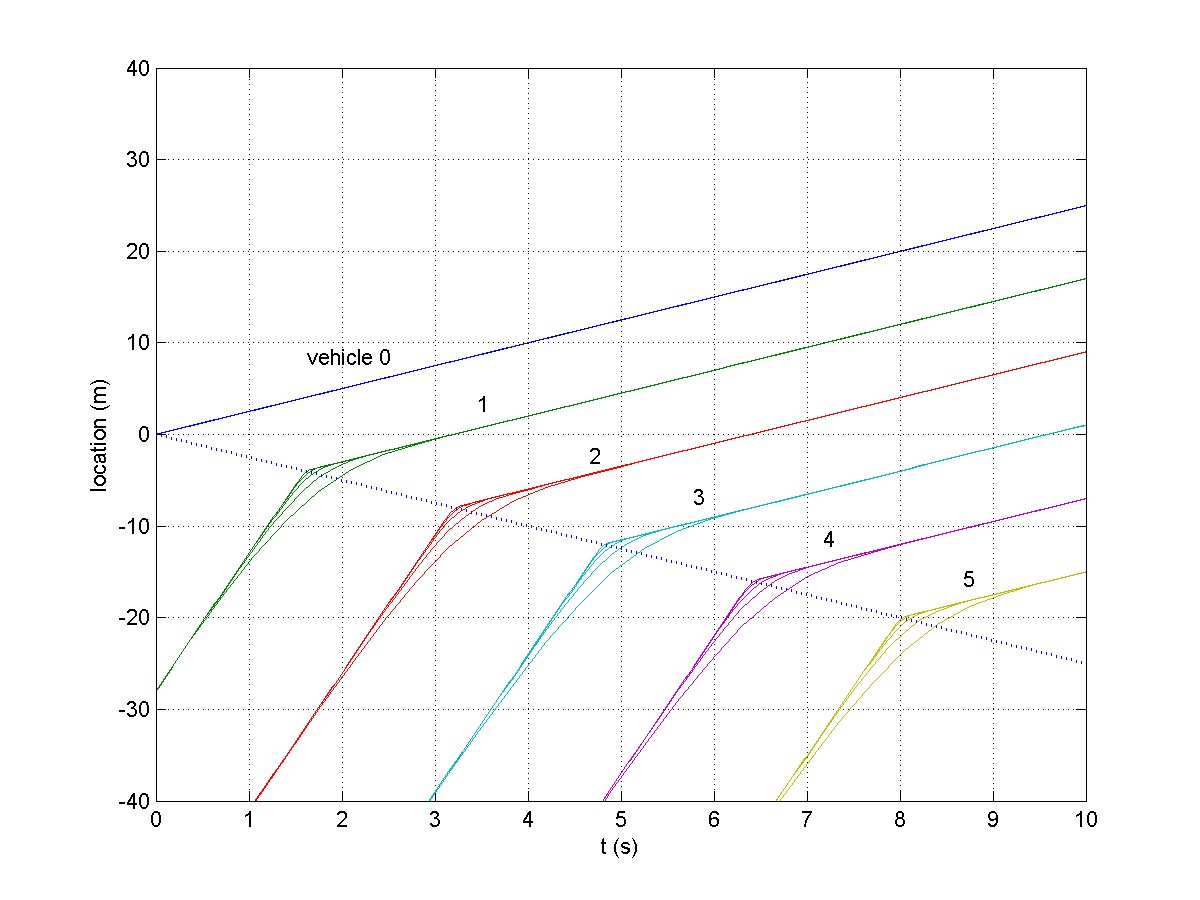} &
\includegraphics[width=3.05in]{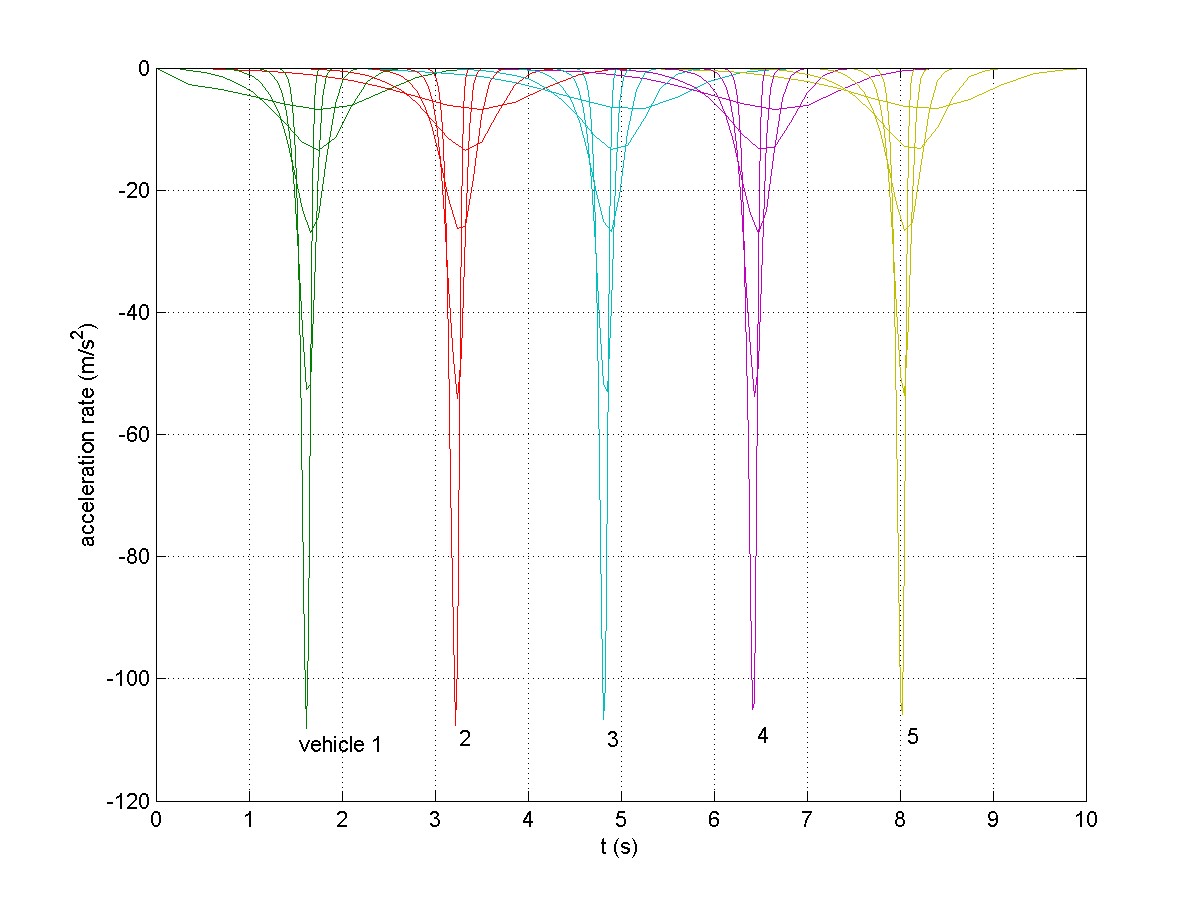} \\
\mbox{\bf (c)} &
\mbox{\bf (d)}
\ea$
\caption{Trajectories and acceleration rates for a platoon of vehicles with $k_1=\frac 14K$: (a) Trajectories when $v_2=\frac 38 V$, (b) Acceleration rates when $v_2=\frac 38 V$, (c) Trajectories when $v_2=\frac 18 V$,  and (d) Acceleration rates when $v_2=\frac 18 V$} \label{lwr2_shock_greenshields} \ec 
\end{figure}

First we set $k_1=\frac 14K$, and $v_1=\frac 34 V$. We consider two cases for the leading vehicle's speed: $v_2=\frac 38 V$ or $\frac 18 V$; in the Riemann problem,  $k_2=\frac 58 K $ or $\frac 78 K$, respectively. Therefore, the LWR model is solved by a shock wave with a speed of $\frac 1 8 V$ or $-\frac 18 V$, correspondingly.
The numerical solutions of the lead-vehicle problem for \refe{discrete_lp2_approx} are shown in \reff{lwr2_shock_greenshields}, where we show five following vehicles' trajectories and acceleration rates along with the leader's, and the dotted line is the theoretical shock wave trajectory. We show trajectories of each vehicle in figures (a) and (c) from bottom to top for $\dN=1,\frac 12,\frac 14,\frac 18$, and $\frac 1{16}$; the acceleration rates are from top to bottom with decreasing $\dN$. From the figures, we have the following observations: (i) the numerical solutions of the trajectories are smooth, compared with theoretical piecewise linear trajectories; but the deceleration rates are very large; (ii)when $\dN$ decreases, all vehicles' trajectories converge, and the deceleration rates become larger; (iii) in the converging trajectories, vehicles  decrease their speeds immediately after the trajectories cross the shock wave, and the deceleration rates  become infinite; (iv) the shock wave speeds are $ \frac 18 V$ and $-\frac 18 V$ in figures (a) and (c), respectively, which are exactly the same as the theoretical predictions; the shock waves can also be seen in the propagation of the deceleration profiles in figures (b) and (d). These results verify that the car-following model indeed converges to the LWR model when $\dt\to 0$, and the deceleration rates are quite large but can be explicitly calculated.

Next we set $k_1=K$ and $v_2=V$. Thus $k_2=0$, and there is no vehicle in front of vehicle 0. This is the queue discharge scenario, and the corresponding Riemann problem is solved by a rarefaction wave with the characteristic wave speed spanning from $V$ to $-V$.
The first five vehicles' trajectories and acceleration rates are shown in \reff{lwr2_qd_greenshields_traj}, where the dotted line is for the characteristic wave with $\phi'(K)=-V$, along which vehicles start to accelerate one by one. We show the trajectories of each vehicle for $\dN=1,\frac 12,\frac 14,\frac 18$, and $\frac 1{16}$ from top to bottom; their acceleration rates are from bottom to top. From the figures we can see that the trajectories get closer with smaller $\dN$, suggesting that the numerical solutions converge. In addition, vehicles' trajectories become smoother along the traffic stream: the first vehicle has the largest acceleration rates (as large as 53.8 $m/s^2$ for $\dN=\frac 1{16}$); and the following vehicles' acceleration rates keep decreasing along the vehicles. 

\begin{figure} \bc
	$\ba{c@{\hspace{0.1in}}c}
	\includegraphics[width=3.2in]{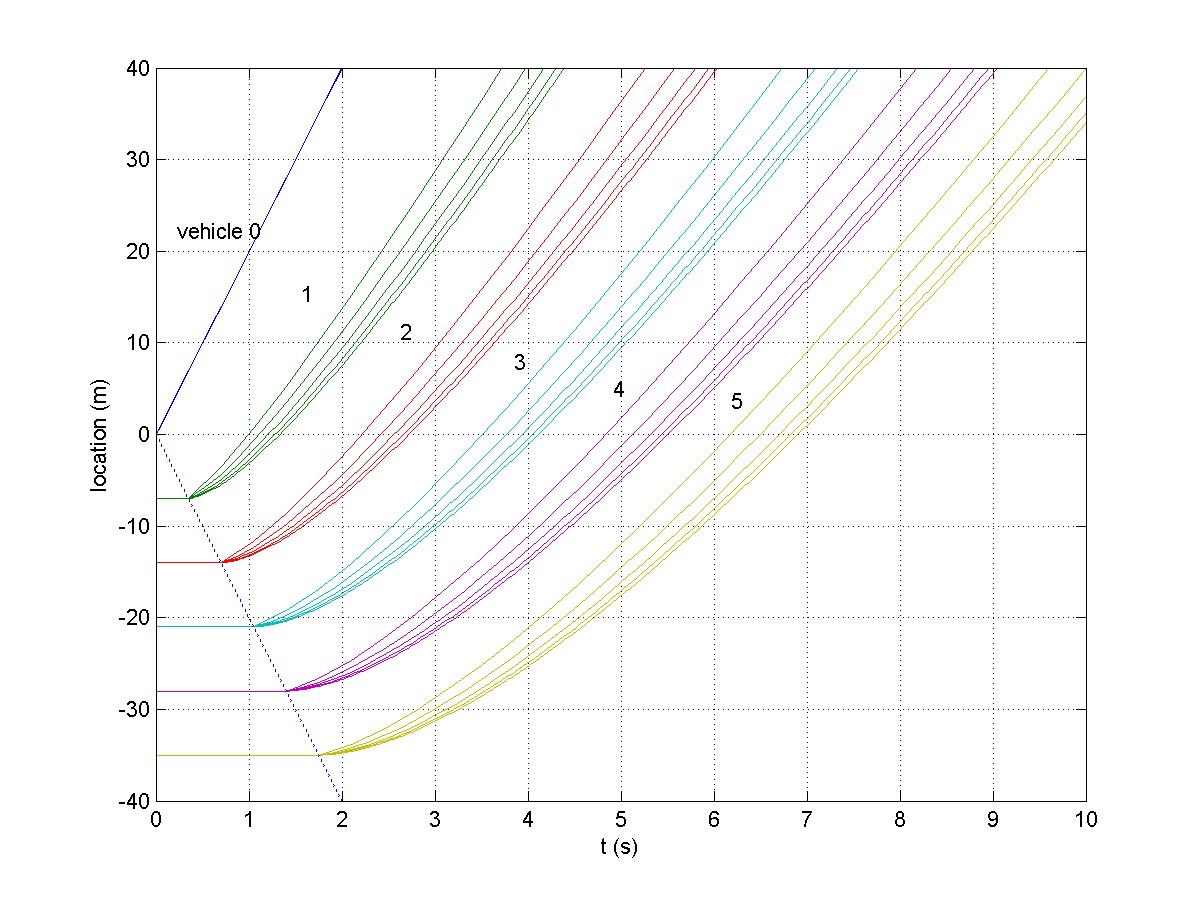} &
	\includegraphics[width=3.05in]{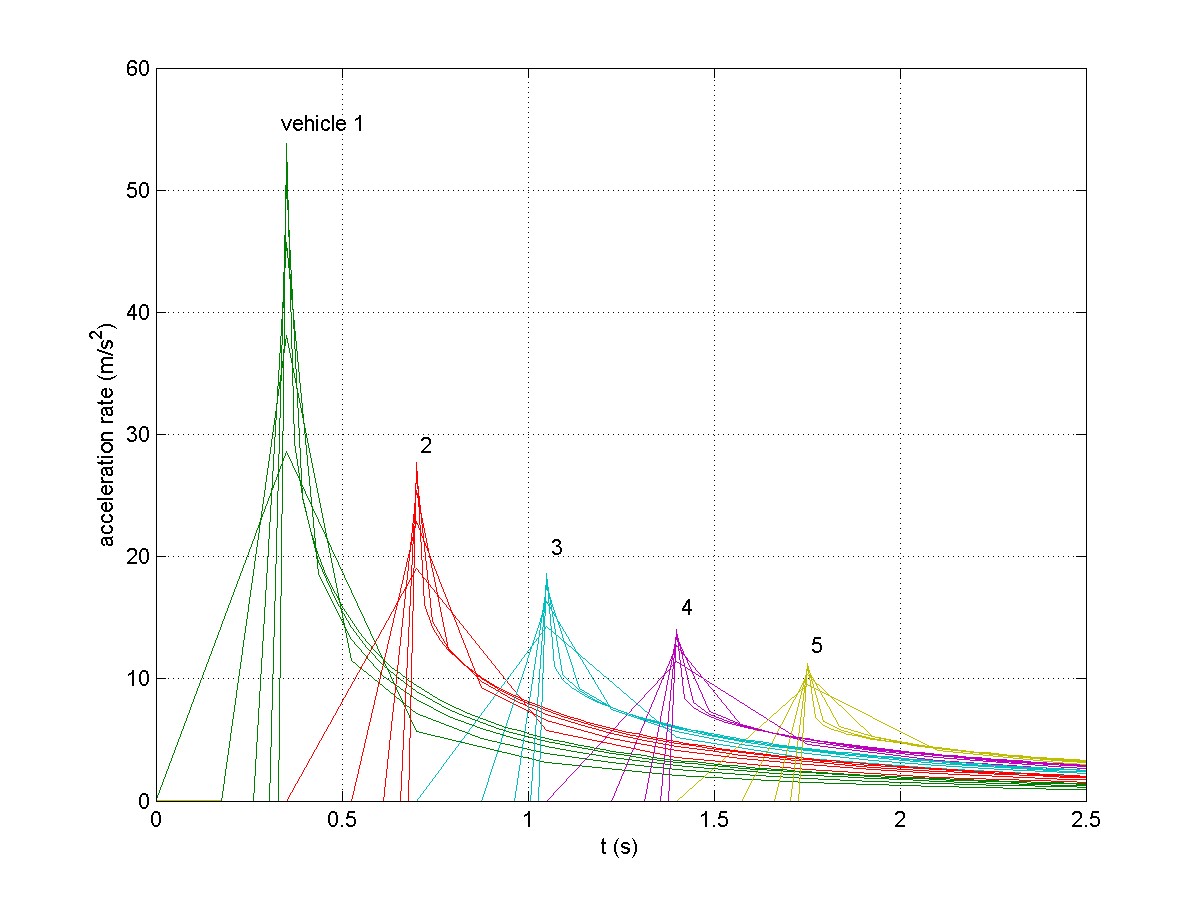} \\
	\mbox{\bf (a)} &
	\mbox{\bf (b)}
	\ea$
\caption{Trajectories and acceleration rates for a platoon of vehicles with $k_1=K$ and $v_2=V$ } \label{lwr2_qd_greenshields_traj}
\ec\efg

\subsection{Triangular fundamental diagram}
In this subsection we solve the lead-vehicle problem with the triangular fundamental diagram, \refe{tri-fd}, where $S=7$ m, $K=1/S$, $V=20$ m/s, and $W=5$ m/s. We set $\dt=1.2 \dN$, which satisfies the CFL condition and the collision-free condition.

First we set $k_1=\frac 1{10}K$, and $v_1= V$. We consider two cases for the leading vehicle's speed: $v_2=\frac 38 V$ or $\frac 1{16} V$. Hence in the Riemann problem $k_2=\frac 25 K $ or $\frac 45 K$, respectively. In this case, the LWR model is solved by a shock wave with a speed of $\frac 16 V$ or $-\frac 1{14} V$, correspondingly.
The numerical solutions are shown in \reff{lwr2_shock_tri}, where we show five following vehicles' trajectories and acceleration rates along with the leader's, and the dotted lines are the theoretical shock wave trajectories. We show trajectories of each vehicle in figures (a) and (c) from bottom to top for $\dN=1,\frac 12,\frac 14,\frac 18$, and $\frac 1{16}$; the acceleration rates are from top to bottom with decreasing $\dN$. From the figures, we have the following observations: (i) in the case of a forward-traveling shock wave, the trajectories are almost the same as the theoretical piecewise linear trajectories for different $\dN$; but the deceleration rates increase with smaller $\dN$; (ii) in the case of a backward traveling shock wave, the trajectories converge to the theoretical piecewise linear trajectories with decreasing $\dN$, and the deceleration rates increase with smaller $\dN$;  (iii) the shock wave speeds are $ \frac 16 V$ and $-\frac 1{14} V$ in figures (a) and (c), respectively, which are exactly the same as the theoretical predictions; the shock waves can also be seen in the propagation of the deceleration profiles in figures (b) and (d).

\begin{figure} \bc
	$\ba{c@{\hspace{0.1in}}c}
	\includegraphics[width=3.1in]{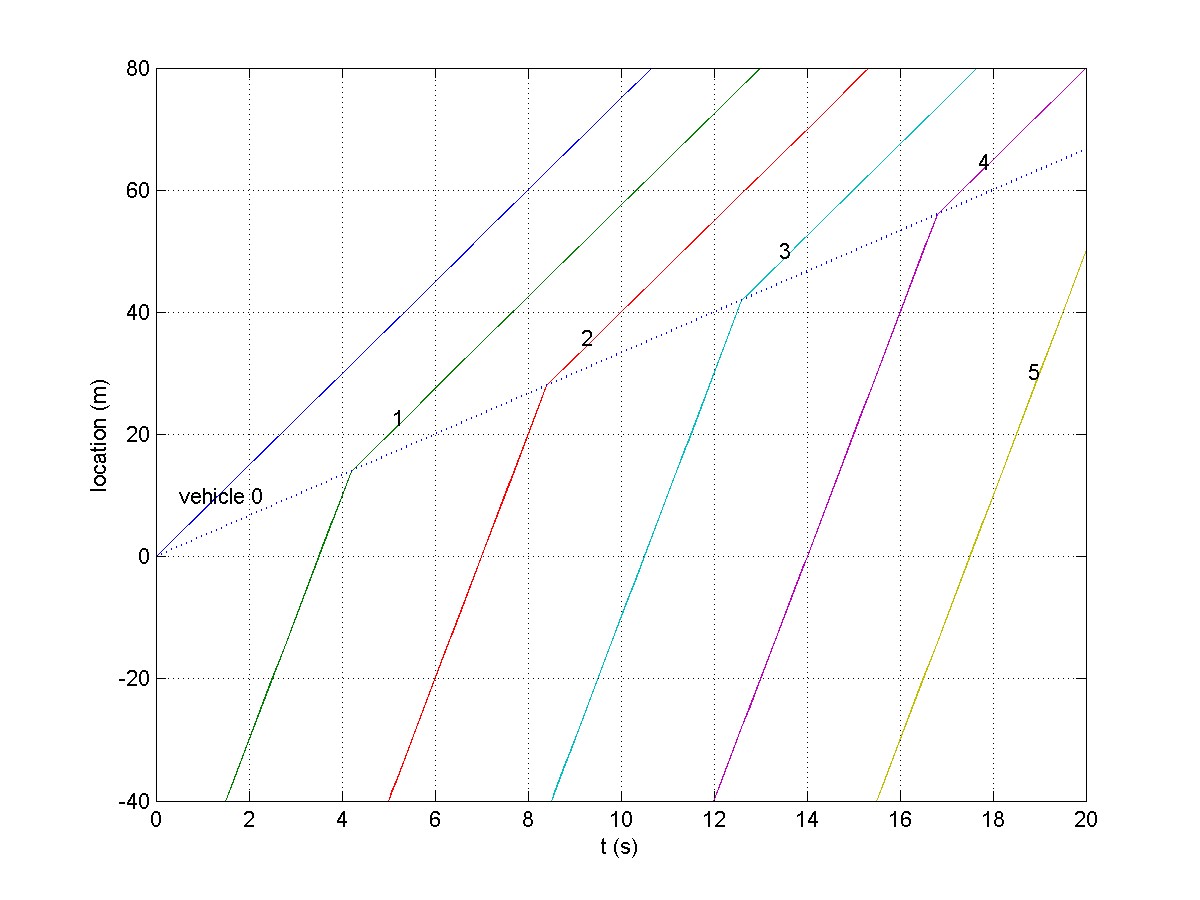} &
	\includegraphics[width=3.1in]{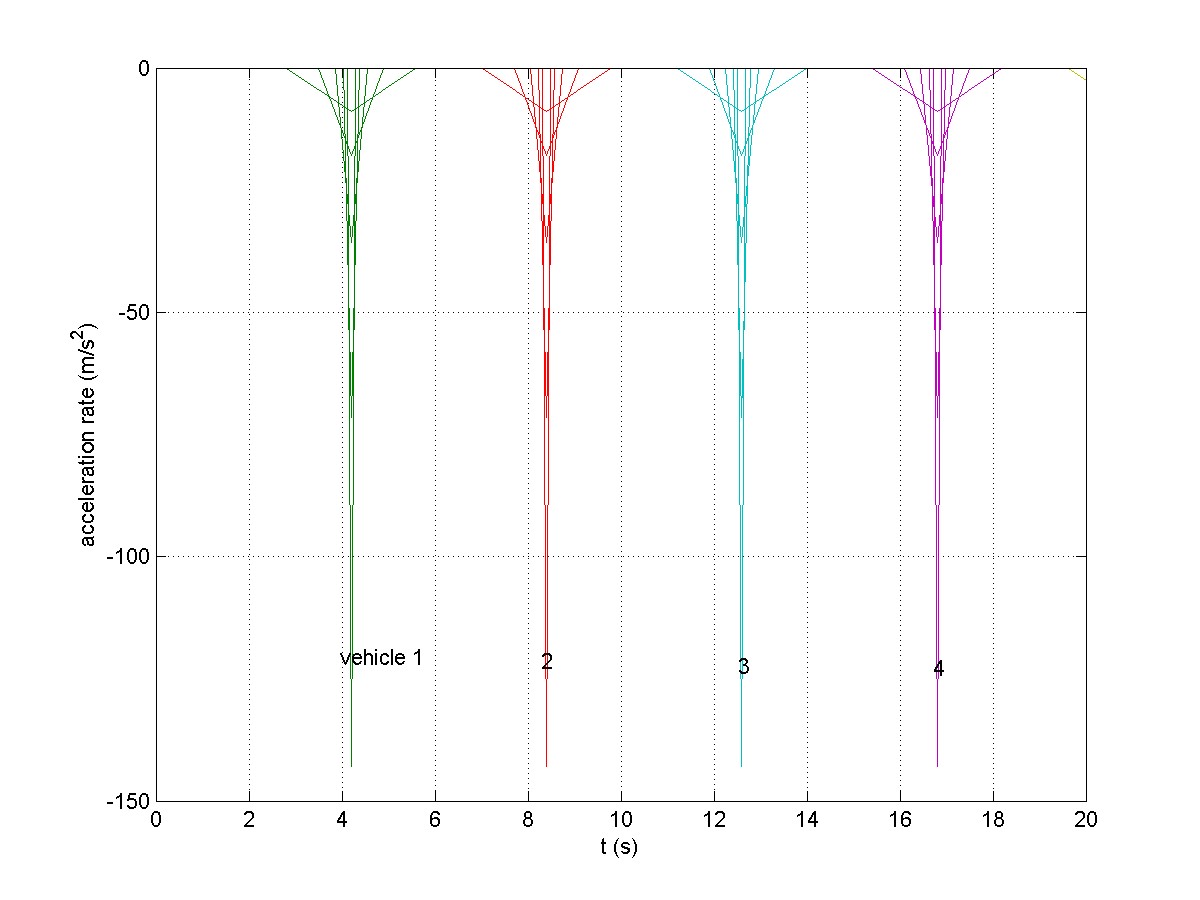} \\
	\mbox{\bf (a)} &
	\mbox{\bf (b)}\\
	\includegraphics[width=3.1in]{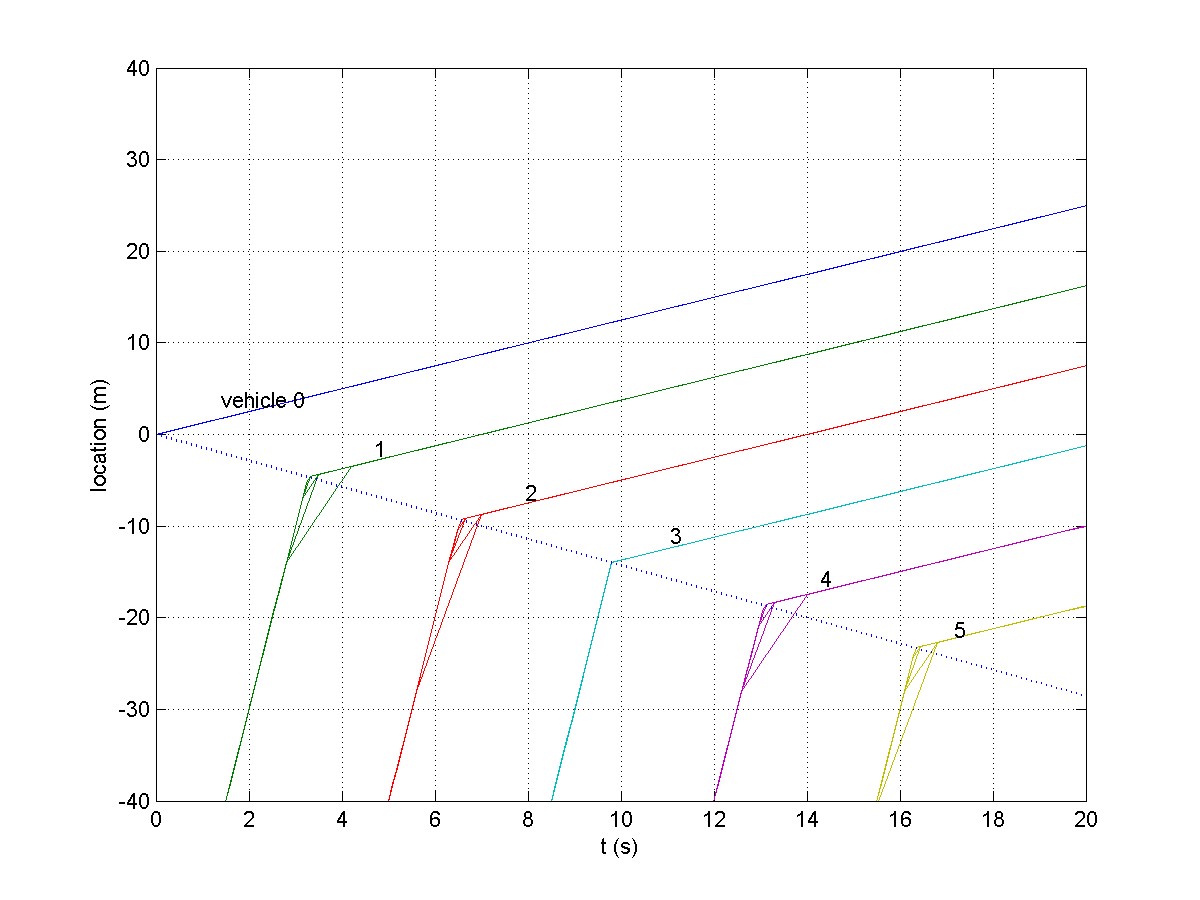} &
	\includegraphics[width=3.1in]{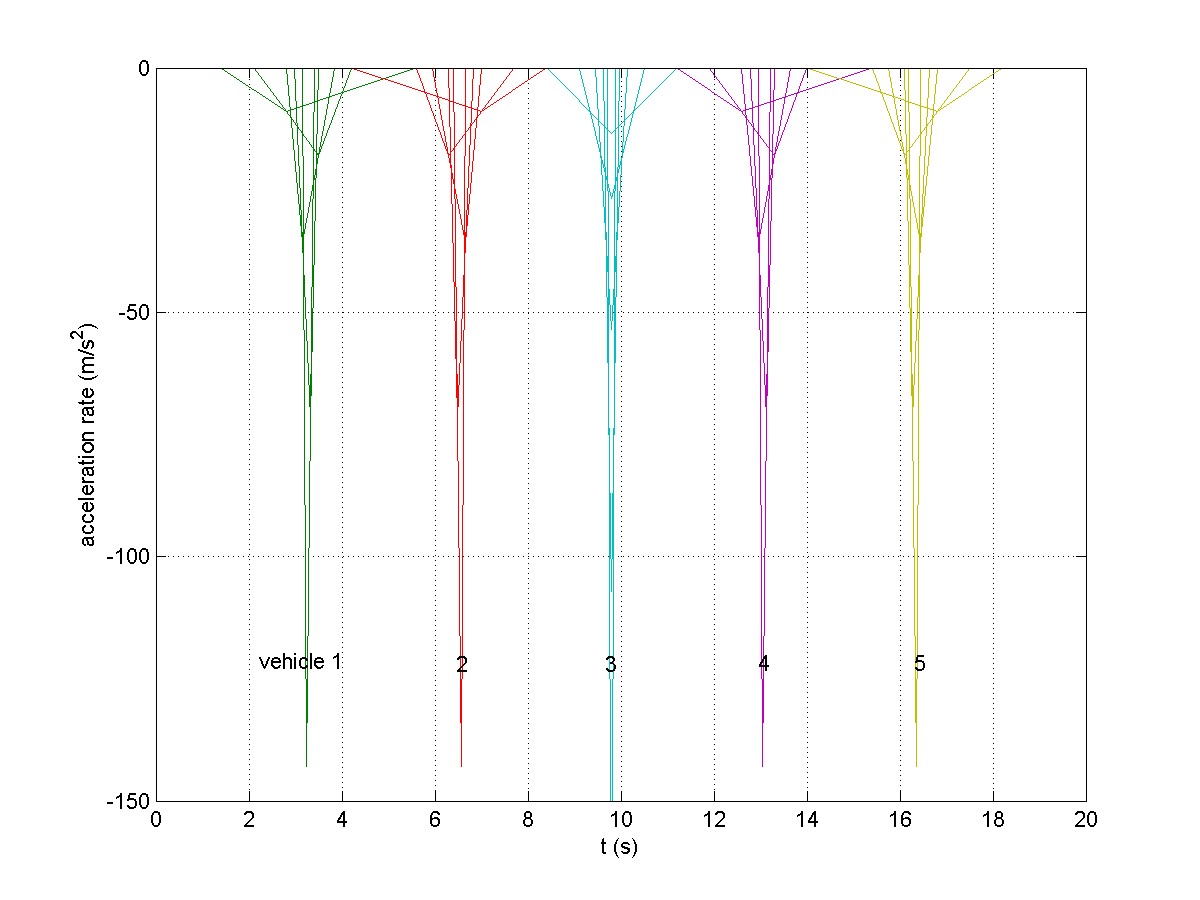} \\
	\mbox{\bf (c)} &
	\mbox{\bf (d)}
	\ea$
	\caption{Trajectories and acceleration rates for a platoon of vehicles with the triangular fundamental diagram and $k_1=\frac 1{10}K$: (a) Trajectories when $v_2=\frac 38 V$, (b) Acceleration rates when $v_2=\frac 38 V$, (c) Trajectories when $v_2=\frac 1{16} V$,  and (d) Acceleration rates when $v_2=\frac 1{16} V$} \label{lwr2_shock_tri} \ec 
\end{figure}

Next we set $k_1=K$ and $v_2=V$. Thus $k_2=0$, and there is no vehicle in front of vehicle 0. This corresponds to the queue discharge scenario, and the corresponding Riemann problem is solved by a degenerate rarefaction wave, which is the same as a shock wave with a speed of $-W$. 
The first five vehicles' trajectories and acceleration rates are shown in \reff{lwr2_qd_tri_traj}, where the dotted line is for the characteristic wave with $\phi'(K)=-W$, along which vehicles start to accelerate one by one. We show trajectories of each vehicle for $\dN=1,\frac 12,\frac 14,\frac 18$, and $\frac 1{16}$; their acceleration rates are from bottom to top. From the figures we can see that the trajectories are almost identical to the analytical solutions with different $\dN$, but the acceleration rates keep increasing with smaller $\dN$, and all vehicles' have the same profile for acceleration rates. Thus the numerical solutions converge to theoretical ones for the LWR model, in which vehicles accelerate instantaneously, and the acceleration rate is infinite.

\begin{figure} \bc
	$\ba{c@{\hspace{0.1in}}c}
	\includegraphics[width=3.1in]{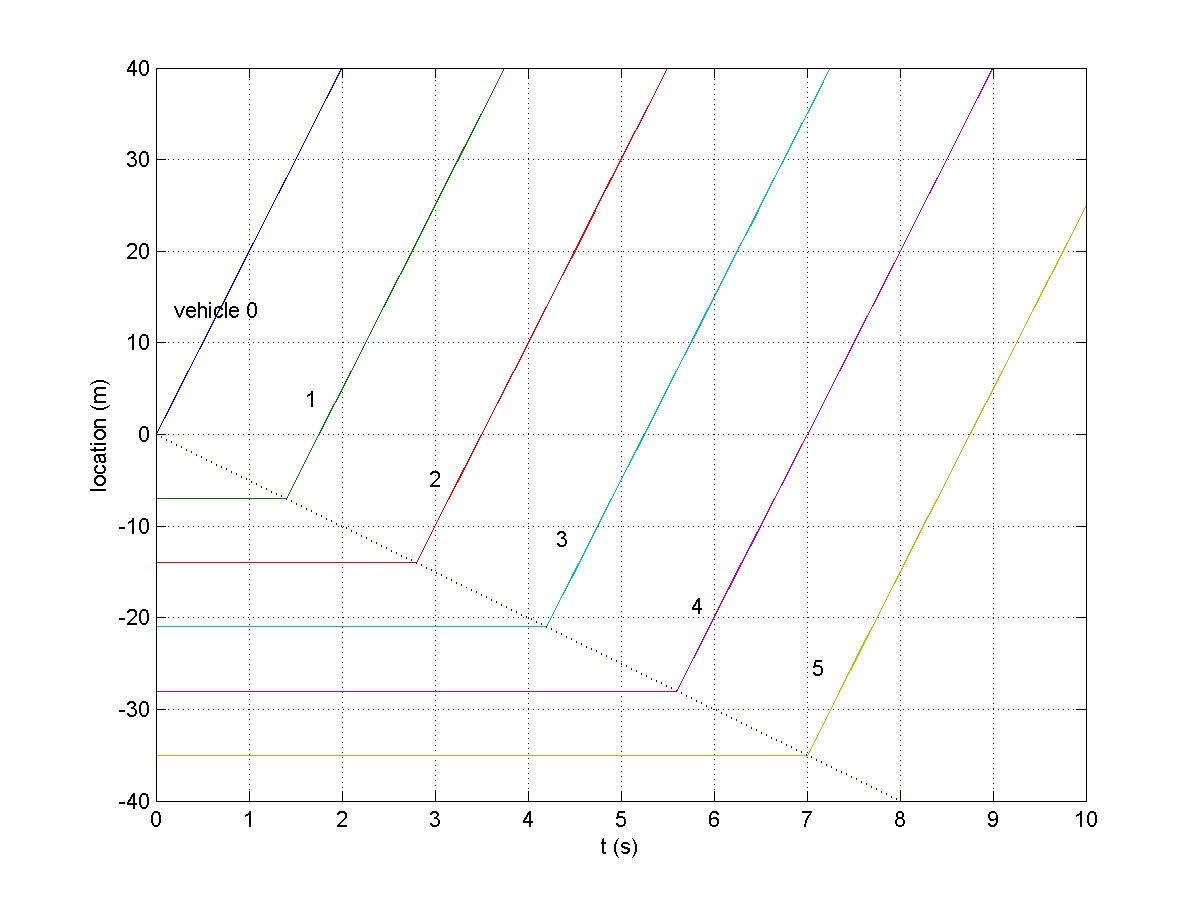} &
	\includegraphics[width=3.1in]{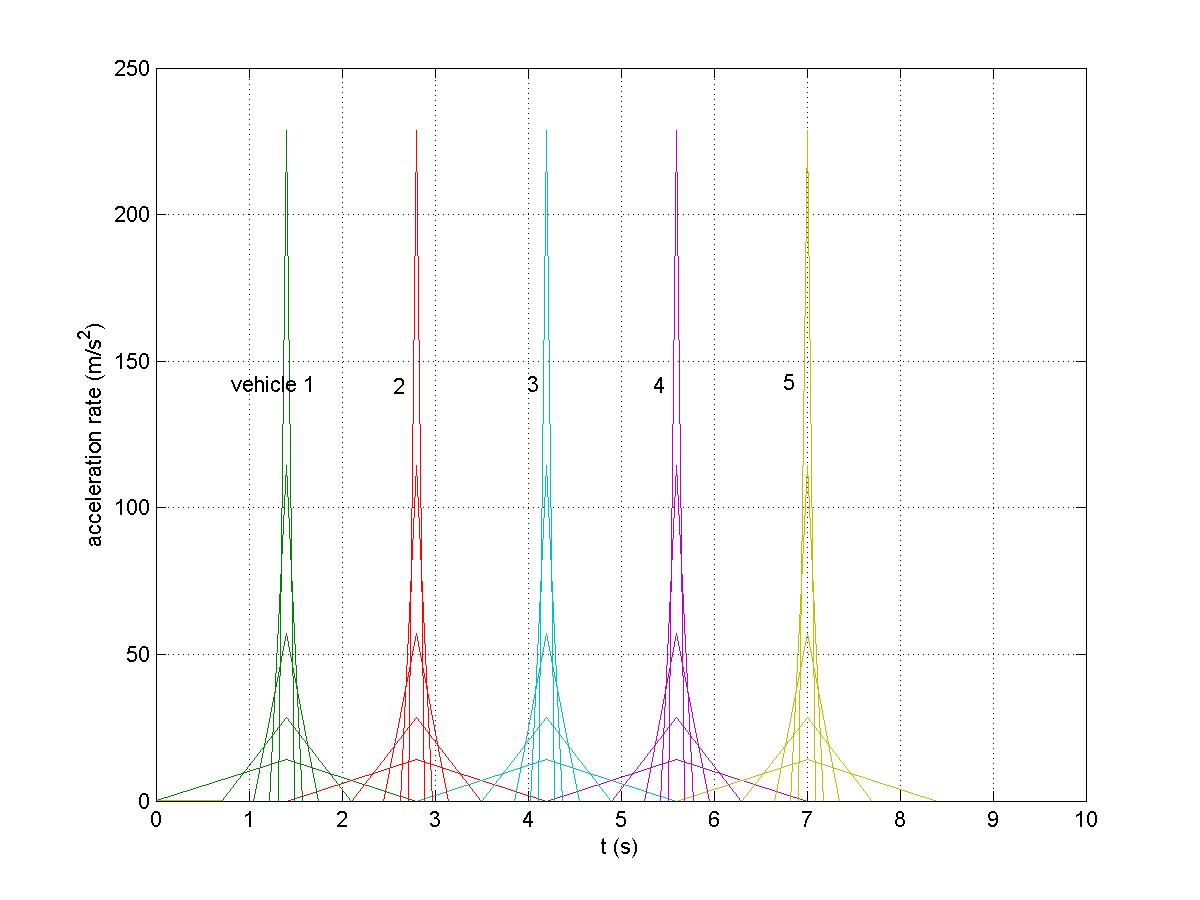} \\
	\mbox{\bf (a)} &
	\mbox{\bf (b)}
	\ea$
	\caption{Trajectories and acceleration rates for a platoon of vehicles with $k_1=K$ and $v_2=V$ } \label{lwr2_qd_tri_traj}
	\ec\efg

\subsection{Non-concave fundamental diagram}
In this subsection, we solve the lead-vehicle problem with a non-concave fundamental diagram, in which the speed-density relation is given by \citep{kerner1994cluster}: $\eta(k)=5.0461[(1+\exp\{[k/K-0.25]/0.06\})^{-1}-3.73\times 10^{-6}] l/T$. Here the unit length $l=28$ m, the relaxation time $T=5$ s, and the jam density $K=0.18$ veh/m. The speed-density relation is decreasing, but the flow-density relation is non-concave. 

From the collision-free condition, \refe{collision-free-cond}, we have $\frac{\dN}\dt \geq 0.89$, or $\dt \leq 1.12 \dN$; in contrast, from the CFL condition, \refe{cfl-cond}, we have $\frac{\dN}\dt \geq 0.32$, or $\dt \leq 3.13 \dN$. Thus the two conditions are not equivalent for such a non-concave fundamental diagram, and the CFL condition is more relaxed.

In the lead-vehicle problem we set $k_1=0.002$ veh/m, and $v_2=0$ m/s. Hence  $k_2=K$ in the corresponding Riemann problem. This corresponds to a scenario in which a very sparse traffic stream with a density of $0.002$ veh/m runs into a red light. We let $\dN=1/10$.

In \reff{lwr2_shock_nonconcave_redlight_2} we show the numerical results for $\dt=\dN$, when both the collision-free and CFL conditions are satisfied, and $\dt=2\dN$, when the CFL condition is satisfied but the collision-free condition is violated. In the top row with $\dt=\dN$, the following vehicles' speeds decrease to zero as in figure (b), and the spacings of all following vehicles decrease to the jam spacing of $1/0.18$ m, and no collisions occur as in figure (a). But in the bottom row with $\dt=2\dN$, when the collision-free condition is violated, vehicles eventually collide into each other as their trajectories overlap in figure (c), and vehicles can develop negative speeds as shown in figure (d). This example confirms the theoretical conclusion that the collision-free condition is necessary for obtaining physically meaningful solutions  in the car-following model, but the CFL condition is no longer sufficient.

\begin{figure} \bc
	$\ba{c@{\hspace{0.1in}}c}
	\includegraphics[width=3.1in]{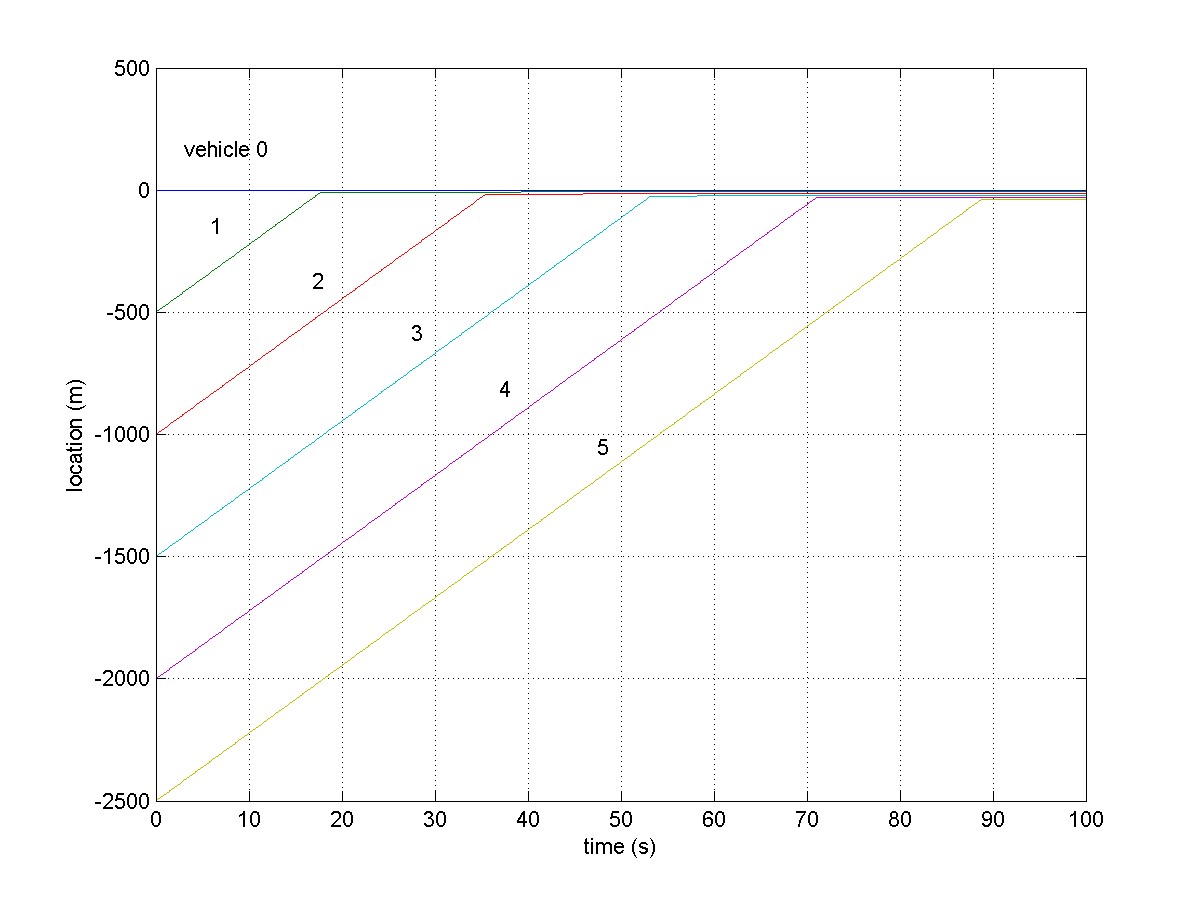} &
	\includegraphics[width=3.1in]{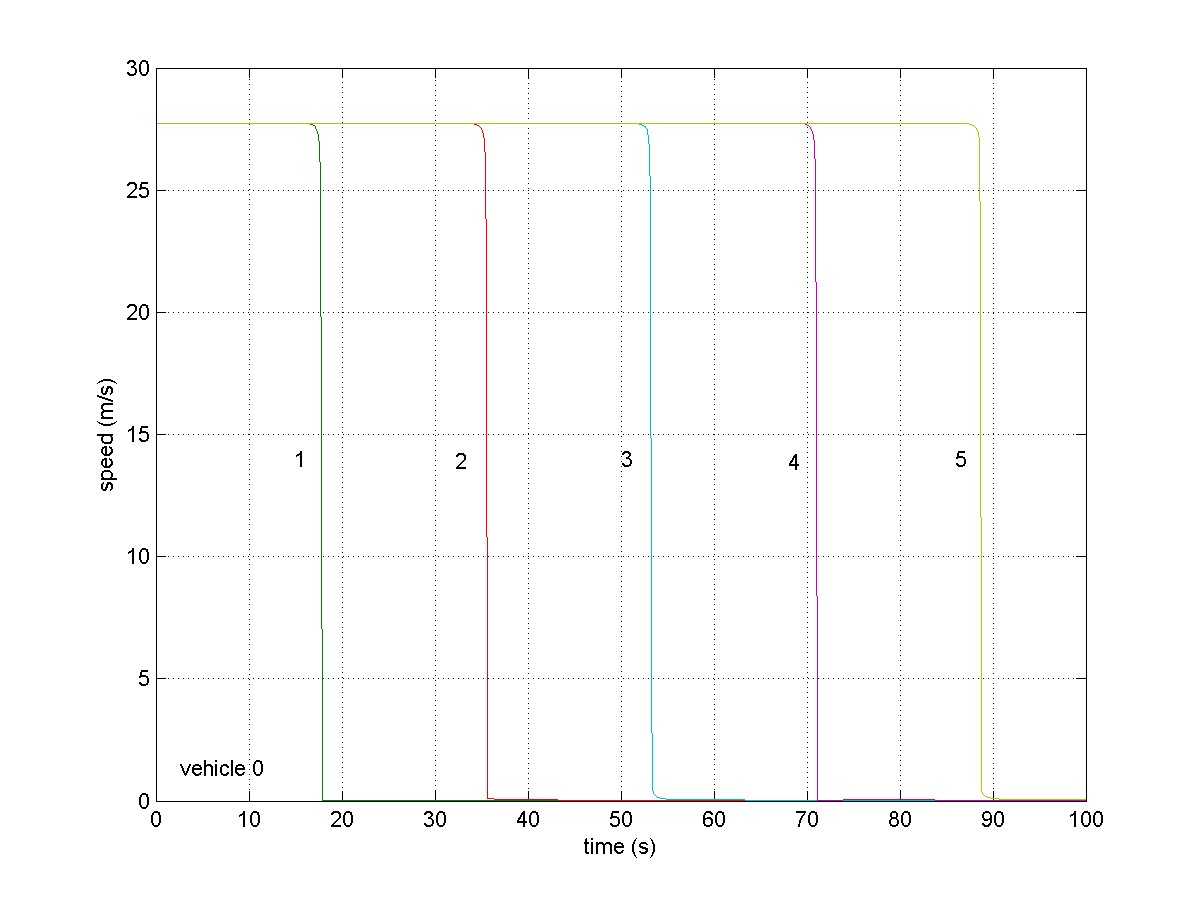} \\
	\mbox{\bf (a)} &
	\mbox{\bf (b)}\\
	\includegraphics[width=3.1in]{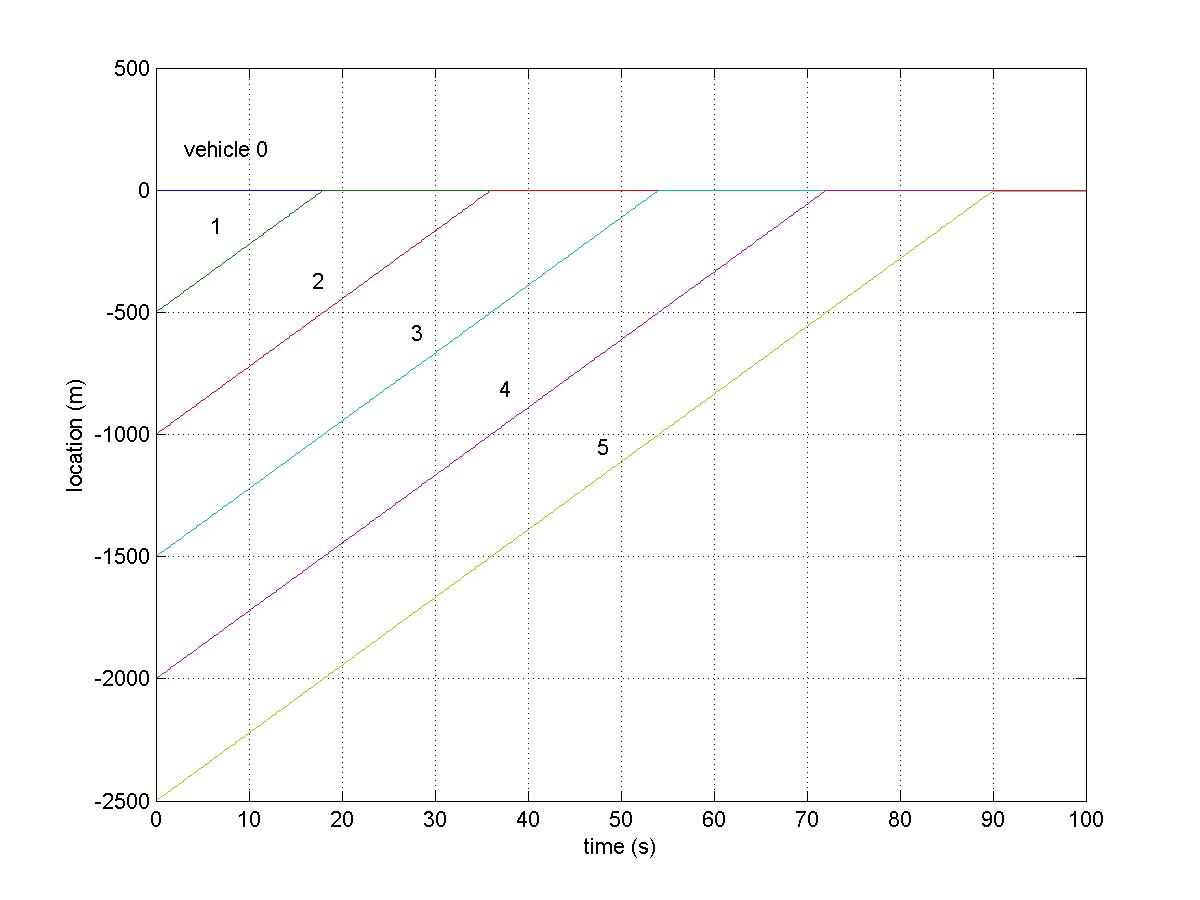} &
	\includegraphics[width=3.1in]{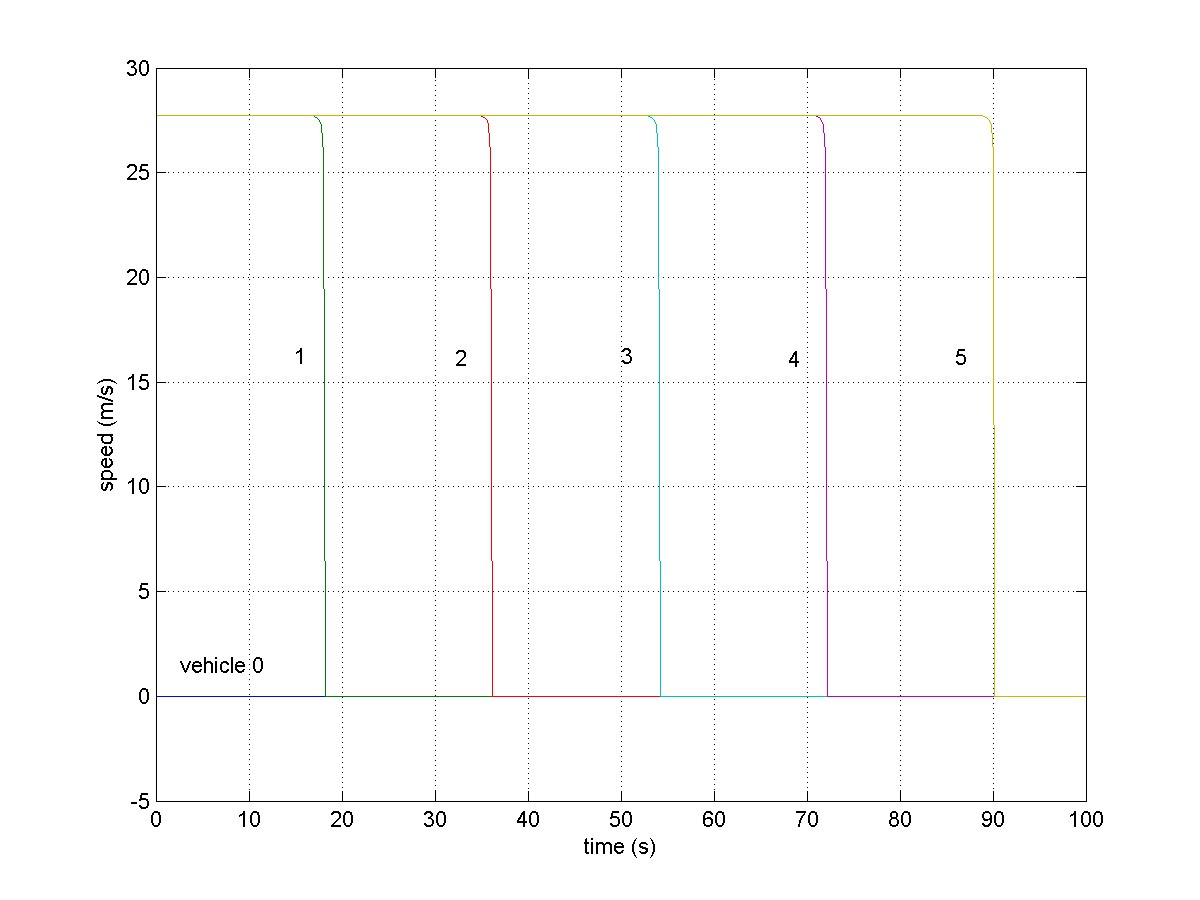} \\
	\mbox{\bf (c)} &
	\mbox{\bf (d)}
	\ea$
	\caption{Trajectories and speeds for the LWR model with a non-concave fundamental diagram: In (a) and (b), $\dt=\dN$ satisfying the collision-free condition; In (c) and (d), $\dt=2\dN$ satisfying the CFL condition, but not the collision-free condition } \label{lwr2_shock_nonconcave_redlight_2}
	\ec\efg

\section{Correction of general second-order models}
In the preceding sections, we established that the time- and vehicle-discrete car-following model, \refe{ns-dis2}, the nonstandard second-order model, \refe{conser_eqn} and \refe{ns-lwr}, and the LWR model, \refe{E-S}, are equivalent and have the following properties, under the collision-free condition given in \refe{collision-free-cond}: (i) existence of a unique speed-density relation: $v=\eta(k)$; (ii) stable; (iii) anisotropic; (iv)  forward-traveling; and (v) collision-free.

In contrast, general second-order models are more realistic than the LWR model, with respect to the instability property and scattered speed-density relations in non-equilibrium traffic. Even though they still admit anisotropic solutions in their equivalent car-following models \citep{jin2016equivalence},  vehicles may travel backwards, as pointed out in \citep{daganzo1995requiem}. An additional undesirable property is that vehicles may collide into each other in these models.
In this section, we present two correction methods  to eliminate such unfavorable properties of general second-order models, but keep their favorable ones. 

Here we consider general second-order continuum models with the following acceleration equation:
\bqn
v_t+vv_x&=&\Psi\left(v,\frac 1k, \frac{v_x}k\right), \label{2ndmodel}
\eqn
in addition to the conservation equation, \refe{conser_eqn}.
We assume that, when in equilibrium; i.e., when the speed is both time- and location-independent, the second-order model has a unique speed-density relation, $v=\eta(k)$. That is, $\Psi\left(\eta(k),\frac 1k,0\right)=0$. Even though this model doesnot have the anticipation term $k_x$ or the viscosity term $v_{xx}$, it is sufficiently general to include continuum models in \citep{Phillips1979traffic,greenberg2001extensions,jiang2002JWZ}  as well as many converted from car-following models by following the method in \citep{jin2016equivalence}.

However, as demonstrated in the following subsection, many general second-order continuum models and their corresponding car-following models may have solutions with negative speeds and collisions.

\subsection{First correction method}
Inspired by the nonstandard second-order model of the LWR model, \refe{conser_eqn} and \refe{ns-lwr}, we present a correction method to adjust \refe{2ndmodel}   as 
\bqn
v_t+vv_x&=&\max\left\{-\frac v \epsilon, \min\left\{\frac{\eta(k)-v}\epsilon, \Psi\left(v,\frac 1k, \frac{v_x}k\right)\right\}\right\}, \label{corrected1}
\eqn
which is equivalent to 
\bqn
v_t+vv_x&=&\min\left\{\frac{\eta(k)-v}\epsilon,\max\left\{-\frac v \epsilon, \Psi\left(v,\frac 1k, \frac{v_x}k\right)\right\}\right\}, \label{corrected2}
\eqn
since $\frac{\eta(k)-v}\epsilon\geq -\frac v \epsilon$. Here we assume that $k\in[0,K]$, and $\eta(k)\geq 0$. We have the following remarks on the correction method:
\ben
\item Both of the correction terms, $-\frac v\epsilon$ and $\frac{\eta(k)-v}\epsilon$ contain the hyperreal infinitesimal number $\epsilon$. Thus the corrected model is a nonstandard second-order model, which may not be equivalent to the LWR model.
\item For the corrected model, the speed-density relation in steady states is still $v=\eta(k)$, the same as that for the original model.
\item Among the two correction terms, $-\frac v\epsilon$ is used to eliminate negative speeds; and $\frac{\eta(k)-v}\epsilon$ is the hyperreal acceleration term in the nonstandard second-order formulation of the LWR model, \refe{ns-lwr}, and used to eliminate collisions. 
\item
In the nonstandard second-order formulation of the LWR model, \refe{conser_eqn} and \refe{ns-lwr},  the acceleration rate $\Psi(v,\frac 1k, \frac{v_x}k)=\frac{\eta(k)-v}\epsilon$, and the corrected model is the same as \refe{conser_eqn} and \refe{ns-lwr}, since $\frac{\eta(k)-v}\epsilon\geq -\frac v \epsilon$. In this sense, the nonstandard second-order formulation, \refe{conser_eqn} and \refe{ns-lwr},  is already corrected.
\een

Following the conversion method in \citep{jin2016equivalence}, the corrected second-order model is equivalent to the following second-order car-following model:
\bsq \label{con-cf-corrected}
\bqn
X_{tt}&=&\max\left\{-\frac {X_t} \epsilon, \min\left\{\frac{\theta(-X_N)-X_t}\epsilon, \Psi(X_t,-X_N, -X_{tN})\right\}\right\},
\eqn
or
\bqn
X_{tt}&=&\min\left\{\frac{\theta(-X_N)-X_t}\epsilon, \max\left\{-\frac {X_t} \epsilon,  \Psi(X_t,-X_N, -X_{tN}) \right\}\right\}.
\eqn
\esq
Further we discretize the model in vehicle $N$ with the anisotropic and symplectic Euler methods and obtain  the time- and vehicle-discrete car-following model:
\bsq \label{cf-corrected1}
\bqn
X_t(t+\dt,N)&=&\max\left\{0,  \min\left\{  \theta\left( \frac{X(t,N-\dN)-X(t,N)}\dN\right), \right.\right.\nonumber\\
&&\left.\left.
 X_t(t,N)+\dt \cdot A(t,N)\right\} \right\}, \label{corrected-speed}\\
 X(t+\dt,N)&=&\max\left\{X(t,N), \min\left\{ X(t,N)+ \dt \cdot \theta\left( \frac{X(t,N-\dN)-X(t,N)}\dN\right), \right.\right. \nonumber\\
 &&\left.\left.
 X(t,N)+ \dt\cdot X_t(t,N)+\dt^2 \cdot A(t,N)\right\} \right\}. \label{corrected-location}
\eqn
\esq
where $A(t,N)$ is the acceleration rate in the original second-order model:
\bqn
A(t,N)=\Psi\left(X_t(t,N), \frac{X(t,N-\dN)-X(t,N)}\dN, \frac{X_t(t,N)-X_t(t,N-\dN)}\dN\right). \label{def:A}
\eqn
Note that the original acceleration rate is determined by vehicle $N$'s speed, the spacing, and the speed difference.

\begin{theorem} When $\dt$ and $\dN$ satisfy the collision-free condition \refe{collision-free-cond}, the corrected car-following models, \refe{corrected-speed} and \refe{corrected-location}, are forward-traveling and collision-free.	
	\end{theorem}
{\em Proof}. From \refe{corrected-speed}, it is straightforward that $X_t(t+\dt,N)\geq 0$; i.e., it is forward-traveling.

From \refe{corrected-location}, we have that
\bqs
X(t+\dt,N-\dN)-X(t+\dt,N)&\geq& X(t,N-\dN)-X(t+\dt,N) \\
&\geq& X(t,N-\dN)-X(t,N) \\&&- \dt \cdot \theta\left( \frac{X(t,N-\dN)-X(t,N)}\dN\right)\\
&=&\frac{\dN}k -\dt \eta(k),
\eqs
where $k=\frac\dN{X(t,N-\dN)-X(t,N)}$.

If \refe{collision-free-cond} is satisfied, then $\frac{\dN}k -\dt \eta(k)\geq S \cdot \dN$, and $X(t+\dt,N-\dN)-X(t+\dt,N)\geq S \cdot \dN$. Thus the model is collision-free.
\eop

\subsection{Second correction method}
For \refe{2ndmodel}, we present another correction method  as 
\bqn
v_t+vv_x&=&\max\left\{-\frac v \epsilon, \min\left\{(\frac 1k -\frac 1 K ) \frac {\delta}{\epsilon^2} -\frac v\epsilon, \Psi\left(v,\frac 1k, \frac{v_x}k\right)\right\}\right\}, 
\eqn
where $\delta=\lim_{\dN\to 0^+} \dN$ is another infinitesimal number and equal to $\dN$ in the discrete form. 

Following \citep{jin2016equivalence}, we obtain the following second-order car-following model:
\bqn
X_{tt}&=&\max\left\{-\frac {X_t} \epsilon, \min\left\{(-X_N-S) \frac {\delta}{\epsilon^2}, \Psi(X_t,-X_N, -X_{tN})\right\}\right\},  \label{con-cf-corrected2}
\eqn
whose discrete form is
\bsq \label{cf-corrected2}
\bqn
X_t(t+\dt,N)&=&\max\left\{0,  \min\left\{  \frac{X(t,N-\dN)-X(t,N) -S\dN }\dt, \right.\right.\nonumber\\
&&\left.\left.
X_t(t,N)+\dt \cdot A(t,N)\right\} \right\}, \label{corrected-speed2}\\
X(t+\dt,N)&=&\max\left\{X(t,N), \min\left\{ X(t,N-\dN) -S\dN, \right.\right. \nonumber\\
&&\left.\left.
X(t,N)+ \dt\cdot X_t(t,N)+\dt^2 \cdot A(t,N)\right\} \right\}. \label{corrected-location2}
\eqn
\esq
where $A(t,N)$ is given in \refe{def:A}. 

We have the following observations regarding the second correction method.
\ben
\item The second correction method was used to study the optimal velocity model in \citep{jin2016equivalence}.
\item \refe{cf-corrected2} is forward-traveling and collision-free for any $\dt$ and $\dN$. 
\item When \refe{collision-free-cond} is satisfied, the speed in the second correction method is greater  than or equal to that in the first correction method, since
\bqs
X(t,N-\dN)-X(t,N) - \dt \cdot \theta\left( \frac{X(t,N-\dN)-X(t,N)}\dN\right) \geq S\dN.
\eqs
\een

\commentout{Is the corrected model still unstable? This correction method is slightly different from that in jin2016equivalence. It seems that the first correction method also stabilizes the model, but not the second one. How to prove this?}

\subsection{An example}
In this subsection, we consider the Jiang-Wu-Zhu (JWZ) second-order model \citep{jiang2002JWZ}, in which the original acceleration equation is given by
\bqn
v_t+vv_x&=&\frac{\eta(k)-v}T+c_0 v_x.
\eqn
Thus its two characteristic wave speeds are $\la_1(k,v)=v$ and $\la_2(k,v)=v-c_0$. If $c_0\leq 0$, it can be shown that the model is always  unstable with the three definitions of stability in Section 2; if $c_0>0$, the condition for it to be stable is $-c_0\leq k\eta'(k) \leq 0$. 
The corresponding anisotropic car-following model is
\bqs
X_t(t+\dt,N)&=& X_t(t,N)+\dt \cdot A(t,N),\\
X(t+\dt,N)&=&X(t,N)+\dt \cdot X_t(t,N)+\dt^2 \cdot A(t,N),
\eqs
where the acceleration rate is given by
\bqs
A(t,N)&=& \frac{ \theta\left( \frac{X(t,N-\dN)-X(t,N)}\dN\right)-X_t(t,N)}T+ c_0  \frac{X_t(t,N-\dN)-X_t(t,N)}{X(t,N-\dN)-X(t,N)) }.
\eqs

Following the first correction method, we obtain the following corrected, nonstandard second-order continuum model:
\bqn
v_t+vv_x&=&\max\left\{-\frac v\epsilon,  \min\left\{\frac{\eta(k)-v}\epsilon , \frac{\eta(k)-v}T+c_0 v_x \right\} \right\},
\eqn
whose car-following formulation is
\bqs
X_t(t+\dt,N)&=&\max\left\{0, \min\left\{  \theta\left( \frac{X(t,N-\dN)-X(t,N)}\dN\right), X_t(t,N)+\dt \cdot A(t,N)\right\} \right\}, \\
X(t+\dt,N)&=&\max\left\{X(t,N), \min\left\{X(t,N)+\dt\cdot  \theta\left( \frac{X(t,N-\dN)-X(t,N)}\dN\right), \right.\right.\\
&&\left.\left. X(t,N)+\dt \cdot X_t(t,N)+\dt^2 \cdot A(t,N)\right\}\right\}.
\eqs

We numerically solve the JWZ model and its corrected version for the following Riemann problem studied in \citep{daganzo1995requiem}:
\bqs
(k(0,x),v(0,x))&=&\cas{{ll} (k_1,0), & x<0\\ (K,0), & x\geq 0}
\eqs
where $k_1$ is the upstream initial density. In Lagrangian coordinates, this is the same the red light scenario studied in Section 5.3. In \citep{daganzo1995requiem}, $k_1=0$; but to  numerically solve the time- and vehicle-discrete version of the JWZ model, we set $k_1$ to be very small, e.g., $1/100 K$, to have vehicles on the upstream part of the road. 
In Lagrangian coordinates, this corresponds to the following lead vehicle problem: $U_m^0=0$ for $m=0,\cdots,M$, and $Y_m^0=-m s_1 \dN$ for $m=0,\cdots,M$, where $s_1=1/k_1$ is the upstream initial spacing. We use the same triangular fundamental diagram as in Section 5.2. In addition, we set $c_0=2$ m/s and $T=5$ s in the JWZ model. In simulations we set $\dt=\dN=1$, which satisfies the CFL condition and, equivalently, the absolute collision-free condition.

The trajectories and speeds of the first five vehicles are shown in \reff{jwz_cf_redlight}, in which the results in the top row are for the original JWZ model, and those in the bottom row for the corrected JWZ model. In the original model, vehicles first accelerate and eventually stop at a jam spacing; but during the process, vehicles can have negative speeds (see figure (b)) and collide into each other (see figure (a)). Thus the original JWZ model is not always forward-traveling or collision-free, even though    it was proved in \citep{jiang2002JWZ} that there is no negative speeds in a special case when $k_1=0$. But both undesirable properties are eliminated in the corrected model (see figures (c) and (d)). Similar results can be observed when we choose other values for $c_0$, including negative ones. These results validates the proposed correction method.

\begin{figure} \bc
	$\ba{c@{\hspace{0.1in}}c}
	\includegraphics[width=3.1in]{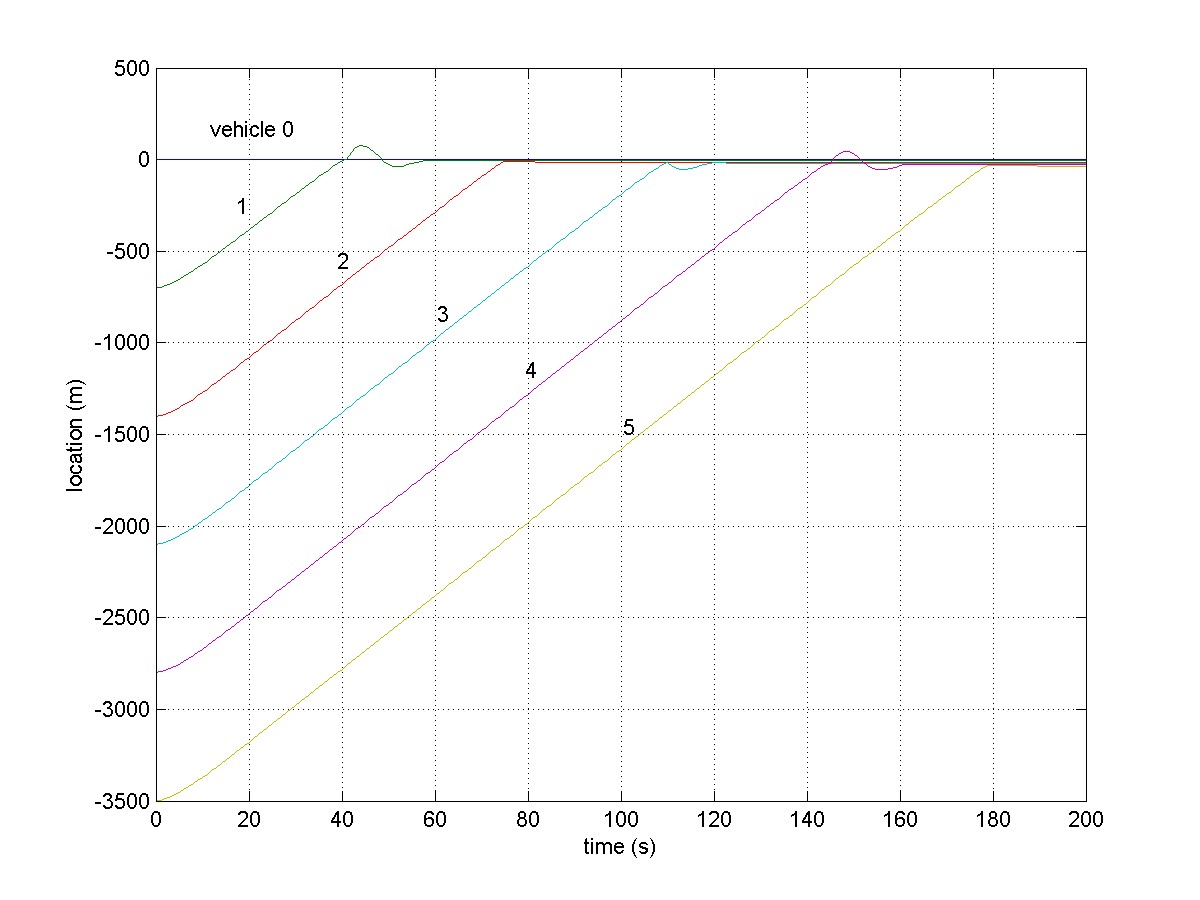} &
	\includegraphics[width=3.1in]{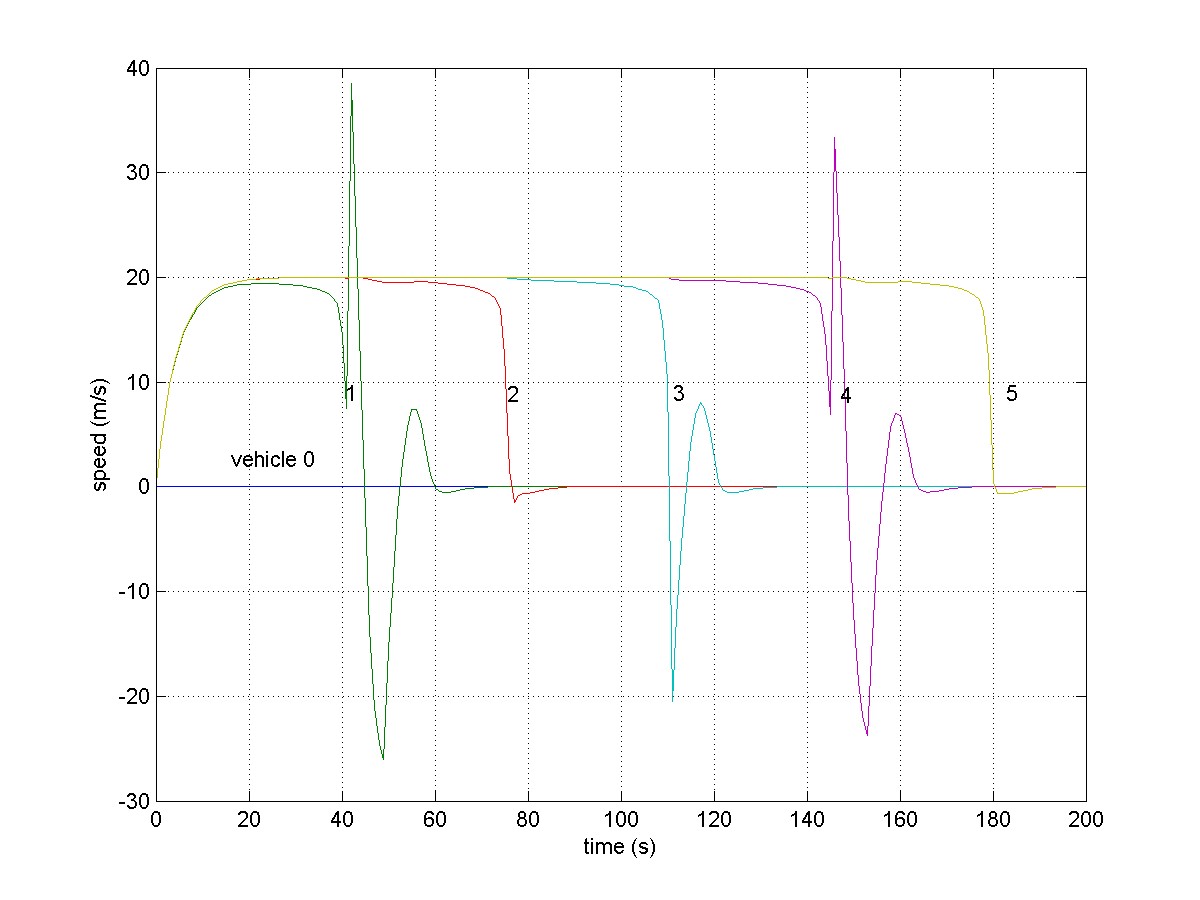} \\
	\mbox{\bf (a)} &
	\mbox{\bf (b)}\\
	\includegraphics[width=3.1in]{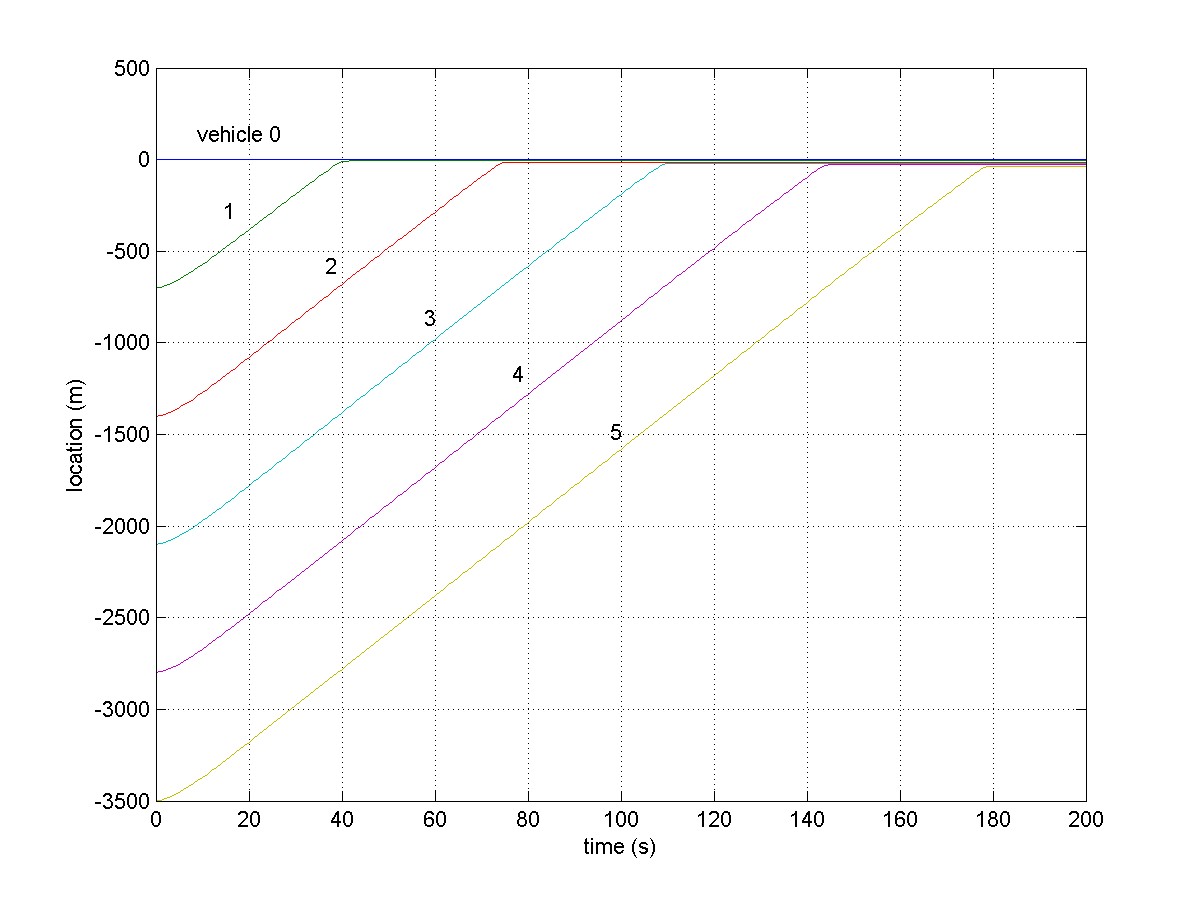} &
	\includegraphics[width=3.1in]{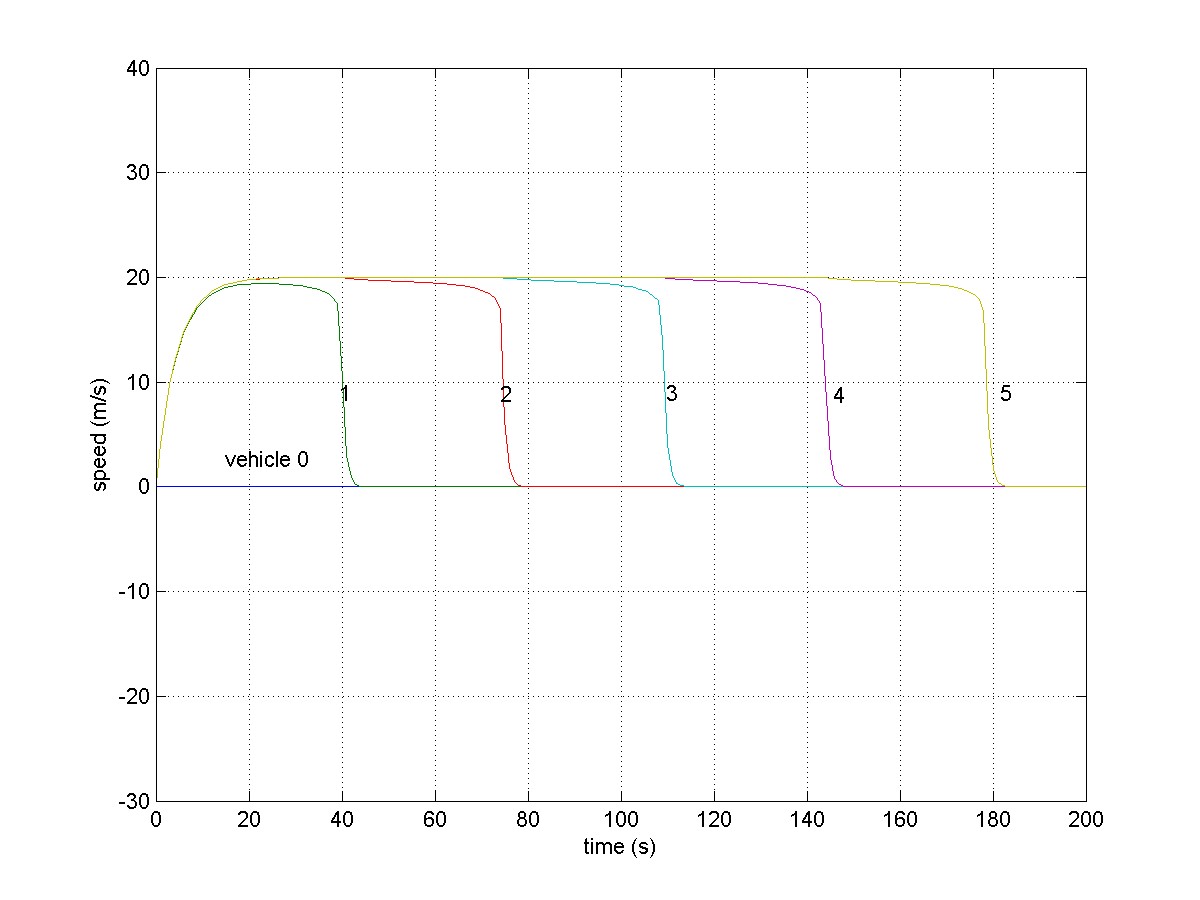} \\
	\mbox{\bf (c)} &
	\mbox{\bf (d)}
	\ea$
	\caption{Trajectories and speeds for the JWZ model: (a) \& (b) before correction; (c) \& (d) after correction method 1 } \label{jwz_cf_redlight}
	\ec\efg

\section{Conclusion}
In this study, we present a nonstandard second-order model by replacing the relaxation time in Phillips' model by a hyperreal infinitesimal number. Since Phillips' model is unstable with three different definitions of stability in both Eulerian and Lagrangian coordinates, we cannot use traditional methods based on the assumption of stability to prove the equivalence between the nonstandard model, which can be considered the zero-relaxation limit of Phillips' model,  and the LWR model,  which is the equilibrium counterpart of Phillips' model. Therefore we resorted to a nonstandard method based on the equivalence relationship between second-order continuum and car-following models in \citep{jin2016equivalence} and proved that the nonstandard model and the LWR model are equivalent since they have the same anisotropic car-following model and stability property. We further derived conditions for the car-following model to be forward-traveling and collision-free, proved that the collision-free condition is consistent with the CFL condition for fundamental diagrams with non-increasing speed-density relations and concave flow-density relations, and demonstrated that other discretization methods, including non-anisotropic discretization methods and the explicit Euler method for both the acceleration and speed, will lead to collisions. With numerical solutions to the lead-vehicle problem we showed that the nonstandard second-order model has the same shock and rarefaction wave solutions as the LWR model for both Greenshields and triangular fundamental diagrams; for a non-concave fundamental diagram we showed that the collision-free condition, but not the CFL condition, yields physically meaningful results. Finally we presented a nonstandard method to correct general second-order models so that their solutions are forward-traveling and collision-free, and verified the method with a numerical example. 

Theoretical and numerical results in this study further validates the equivalence relationship between higher-order continuum and car-following models established in \citep{jin2016equivalence}. For hyperbolic conservation laws with non-convex fluxes,  non-classical shock wave theories have been developed \citep{hayes1997nonclassical,lefloch2002hyperbolic}. In the same spirit, Oleinik and Lax's entropy conditions should also be revised  for  the LWR model with non-concave fundamental diagrams or second-order continuum models. In particular, physically meaningful solutions of a continuum model should be consistent with those of the equivalent car-following model. 
In a sense, the  car-following model can serve as an entropy condition to pick out unique and physical solutions for the corresponding continuum model.

In \citep{jin2016equivalence}, we demonstrated that all second-order continuum models admit anisotropic solutions, through their equivalent car-following models. Therefore, anisotropy is not an intrinsic property of a continuum model, rather it is determined by the discretization method. In this study, we verified that Daganzo's second critique regarding the existence of negative speeds is valid for second-order continuum and car-following models. Moreover, some second-order models are not collision-free. Both undesirable properties were demonstrated for a second-order continuum model and the corresponding car-following model in Section 6.2. However, we showed that all second-order continuum and car-following models can be corrected with the method developed in Section 6.1, so as to eliminate negative speeds and collisions. Therefore, \citep{jin2016equivalence} and this study present a new approach to address the two critiques on second-order continuum models in \citep{daganzo1995requiem}: we can follow the correction method and convert them into car-following models to obtain anisotropic solutions that are forward-traveling and collision-free. 

In this study, we demonstrated that the nonstandard second-order formulation of the LWR model and, equivalently, the car-following model,  have a unique speed-density relation in steady states, is stable, and has anisotropic, forward-traveling, and collision-free solutions. This set of physical properties can serve as conceptual guidelines for analyzing, discretizing, and developing higher-order continuum and car-following models. For example, physically meaningful models and discretization methods should lead to a unique speed-density relation in steady states as well as forward-traveling and collision-free solutions. In particular, the collision-free condition is critical for choosing different discretization methods as well as the time- and vehicle-step sizes; we also demonstrated that it is more general than the CFL condition, which was necessary for converging numerical methods for hyperbolic conservation laws. In particular, such a collision-free condition is critical  for developing car-following algorithms to control connected and automated vehicles. However, the list of physical properties is not complete. In the future we will be interested in introducing upper and lower bounds of the acceleration rate based on the nonstandard second-order formulation of the LWR model as well as the equivalence relationship between continuum and car-following models.

\end{document}